\documentclass{elsart}
\usepackage{amssymb}
\usepackage{amsmath}
\usepackage[english]{babel}
\newtheorem{theorem}{Theorem}

\newtheorem{corollary}[theorem]{Corollary}
\newtheorem{definition}[theorem]{Definition}

\newtheorem{lemma}[theorem]{Lemma}

\newtheorem{proposition}[theorem]{Proposition}
\newtheorem{remark}[theorem]{Remark}
\begin{document}
\begin{frontmatter}
\title{A semi-classical trace formula at a non-degenerate critical level.}
\author{Brice Camus}
\thanks[label2]{This work was supported by the
\textit{TMR-Network PDE \& QM}, reference number ERBFMRXCT960001
and \textit{IHP-Network}, reference number HPRN-CT-2002-00277.}
\address{Mathematisches Institut der Ludwig-Maximilians-Universit\"{a}t M\"{u}nchen.\newline
Theresienstra\ss{}e 39, 80333 Munich, Germany.\newline Email :
brice.camus@univ-reims.fr}
\begin{abstract}
We study the semi-classical trace formula at a critical energy
level for a $h$-pseudo-differential operator whose principal
symbol has a unique non-degenerate critical point for that energy.
This leads to the study of Hamiltonian systems near equilibrium
and near the non-zero periods of the linearized flow. The
contributions of these periods to the trace formula are expressed
in terms of degenerate oscillatory integrals. The new results
obtained are formulated in terms of the geometry of the energy
surface and the classical dynamics on this surface.
\end{abstract}
\end{frontmatter}
\mathindent=0pt
\section{Introduction}
Let be $P_{h}$ a $h$-pseudodifferential, or more generally
$h$-admissible (see \cite{Rob}), self-adjoint operator on
$\mathbb{R}^{n}$. The semi-classical trace formula studies the
asymptotic behavior, as $h$ tends to 0, of the sums
\begin{equation}
\gamma (E,h,\varphi )=\sum\limits_{|\lambda _{j}(h)-E|\leq
\varepsilon }\varphi (\frac{\lambda _{j}(h)-E}{h}), \label{Def
trace}
\end{equation}
where the $\lambda _{j}(h)$ are the eigenvalues of $P_{h}.$ Here
we suppose that the spectrum is discrete in $[E-\varepsilon
,E+\varepsilon ]$, some sufficient conditions for this will be
given below. Let be $p$ the principal symbol of $P_{h}$ and
$\Phi_{t}$ the Hamilton flow of $p$. The semi-classical trace
formula establishes a link between the asymptotic behavior of
(\ref{Def trace}), as $h\rightarrow 0$, and the closed
trajectories of $\Phi _{t}$ of energy $E$. An energy $E$ is said
to be regular when $\nabla p(x,\xi )\neq 0$ on $\Sigma _{E}$,
where $\Sigma _{E}=\{(x,\xi )\text{ }/\text{ }p(x,\xi )=E\}$ is
the surface of energy\ level $E$, and critical if it is not
regular. The case of a regular energy has been intensively studied
and explicit expressions in term of $\Phi_{t}$ are known for the
leading term of (\ref{Def trace}), under suitable conditions on
the flow and when the Fourier transform $\hat{\varphi}$ of
$\varphi$ is supported near a period, see, e.g., Gutzwiller
\cite{GUT}, Balian and Bloch \cite{BB} for the physical
literature, and from a mathematical point of view Brummelhuis and
Uribe \cite{BU}, Petkov and Popov \cite{P-P}, Charbonnel and Popov
\cite{C-P}, Paul and Uribe \cite{PU}.

Here we are interested in the case of a critical energy $E_{c}$ of
$p.$ Brummelhuis, Paul and Uribe in \cite{BPU} have studied the
semi-classical trace formula at a critical energy for quite
general operators but limited to ''small times'', that is for
$\rm{supp}(\hat{\varphi})$ contained in such a small neighborhood
of the origin that the only period of the linearized flow in
$\rm{supp}(\hat{\varphi})$ is 0. Khuat-Duy, in \cite{KhD},
\cite{KhD1}, has obtained the contributions of the non-zero
periods of the linearized flow for arbitrary $\varphi$ with
compactly supported $\hat{\varphi}$, for Schr\"{o}dinger operators
$-\Delta +V(x)$ with $V(x)$ a non-degenerate potential. In this
case the main contribution of such a period was obtained as a
regularization of the Duistermaat-Guillemin density $\rho
_{t}(x,\xi)=|\det (d\Phi _{t}(x,\xi)-\mathrm{Id}|^{-\frac{1}{2}}$.
Generalizing Khuat-Duy's result to more general operators was an
open problem and is the purpose of this article. For a critical
energy level $E_c$ of an arbitrary $h$-admissible operator, and
$(x_{0},\xi _{0})$ a critical point of $p$, the closed
trajectories of the linearized flow
\begin{equation}
\phi _{t}(u)=d_{x,\xi }\Phi _{t}(x_{0},\xi _{0})u=u,\text{ }u\in
T_{(x_{0},\xi _{0})}(T^{\ast }\mathbb{R}^{n}),\text{ }t\neq 0,
\label{sing}
\end{equation}
are not necessarily generated by a positive quadratic form,
contrary to the case of a Schr\"{o}dinger operator, and will give
rise to new contributions to the semi-classical trace formula,
other than those obtained in \cite{BPU}, \cite{KhD} and
\cite{KhD1}. More precisely, viewing $d_{x,\xi }\Phi_{t}(x_{0},\xi
_{0})$ as the Hamiltonian flow of the Hessian $d^{2}p(x_{0},\xi
_{0}),$ here interpreted as an intrinsic quadratic form, these new
contributions to the trace formula arise from the non-trivial
closed trajectories of $d_{x,\xi }\Phi_{t}(x_{0},\xi _{0})$ of
zero energy
\begin{equation}
\left\{
\begin{matrix}
\phi _{T}(u)=u,\text{ }T\neq 0,\\
d^{2}p(x_{0},\xi _{0})(u)=0,\text{ } u\neq 0.
\end{matrix}
\right. \label{E=0}
\end{equation}
The reader can easily verify that the set of such $u$ is empty in
the case of a Schr\"{o}dinger operator : this explains why the new
contributions obtained here does not appear for these operators.
We will show that the new contributions are supported in
$(E_{c},x_{0},\xi _{0},T,u)$ with $(T,u)$ satisfying (\ref{E=0})
and can be expressed in term of $d^{2}p(x_{0},\xi _{0})$ and
higher order derivatives of the flow in $(x_{0},\xi _{0})$.
\section{General hypotheses and main results}
Let $P_{h}$ be in the class of $h$-admissible operators on
$\mathbb{R}^{n}$ with a real symbol. We refer to \cite{Rob} for
the principal notions of semi-classical analysis which we will
use. We note $p$ the principal symbol of $P_{h},$ and $p^{1}$ the
sub-principal symbol. For $E_{c}$ a critical energy of $p$ we will
study the asymptotic behavior of the spectral function $\gamma
(E_{c},h)$, defined by
\begin{equation}
\gamma (E_{c},h)=\sum\limits_{\lambda _{j}(h)\in \lbrack
E_{c}-\varepsilon ,E_{c}+\varepsilon ]}\varphi (\frac{\lambda
_{j}(h)-E_{c}}{h}),
\end{equation}
under the following hypotheses which are classical in this context
:
\medskip
\newline $(H_{1})$ $\exists \varepsilon_{0}>0$\textit{\ such that
}$p^{-1}([E_c-\varepsilon_{0},E_c+\varepsilon_{0}])$\textit{\ is
compact.}
\newline $(H_{2})$\textit{\ }$z_{0}=(x_{0},\xi
_{0})$\textit{\ is the unique critical point of }$p$ \textit{on
the energy surface }$\Sigma _{E_{C}}$\textit{.}
\newline $(H_{3})$ $z_{0}$\textit{\ is a non degenerate critical
point of $p$ and its Hessian $d^{2} p(z_{0})$ is diagonal in some
suitable set of local symplectic coordinates near }$z_{0}$
\begin{equation}\label{formloc}
p(x,\xi)=E_{c}+\frac{1}{2}\sum_{j=1}^{n}w_{j}((x_{j}-x_{0,j})^{2}+\sigma
_{j}(\xi _{j}-\xi _{0,j})^{2})+\mathcal{O}(||(x-x_{0},\xi -\xi
_{0})||^{3}),
\end{equation}
\textit{with }$\sigma _{j}=\pm 1$\textit{\ and }$w_{j}\in \mathbb{R}%
\backslash \{0\}.$ \newline
$(H_{4})$\textit{\ }$\varphi $ \textit{is in the Schwartz space }$\mathcal{S}%
(\mathbb{R})$\textit{\ with Fourier transform }$\hat{\varphi}\in
C_{0}^{\infty }(\mathbb{R})$.

By a classical result, see e.g. \cite{Rob}, the hypothesis
$(H_{1})$ insures that the spectrum of $P_{h}$ is discrete in
$I_{\varepsilon }=[E_{c}-\varepsilon ,E_{c}+\varepsilon ]$ for
$\varepsilon<\varepsilon_{0}$ and $h$ small enough : this will be
assumed in the following. We note $\mathrm{Exp}(tH_{f}),$ with
$H_{f}=\partial _{\xi }f.\partial_{x} -\partial
_{x}f.\partial_{\xi},$ the Hamilton flow of a function $f\in
C^{\infty }(T^{\ast }\mathbb{R}^{n})$. For $\Phi
_{t}=\mathrm{Exp}(tH_{p})$ taking the derivative with respect to
the initial conditions gives a symplectomorphism $d_{x,\xi }\Phi
_{t}(x,\xi ) : T_{(x,\xi )}(T^{\ast }\mathbb{R}^{n})\rightarrow
T_{\Phi _{t}(x,\xi )}(T^{\ast }\mathbb{R}^{n})$, and for
$z_{0}=(x_{0},\xi _{0})$ a critical point of $p$ we have the
fundamental automorphism
\begin{equation}
d_{x,\xi }\Phi _{t}(z_{0}):T_{z_{0}}(T^{\ast
}\mathbb{R}^{n})\rightarrow T_{z_{0}}(T^{\ast }\mathbb{R}^{n}).
\end{equation}
Near the critical point of $p$ we can write
\begin{equation*}
p(x,\xi )=E_{c}+\sum\limits_{j=2}^{N}p_{j}(x,\xi
)+\mathcal{O}(||(x-x_{0},\xi -\xi _{0})||^{N+1}),
\end{equation*}
where the functions $p_{j}$ are homogeneous of degree $j$ in
$(x-x_{0},\xi -\xi _{0}).$ In particular, $p_{2}$ is the Hessian
in $z_{0}$ and can be interpreted as an invariantly defined
quadratic form on $T_{z_{0}}(T^{\ast }\mathbb{R}^{n})$.
\begin{definition}
For all $T\in \mathbb{R}$ let be $\frak{F}_{T}=
\mathrm{Ker}(d_{x,\xi }\Phi _{T}(z_{0})-\mathrm{Id}) \subset
T_{z_{0}}(T^{\ast }\mathbb{R}^{n})$ and
$\mathbb{\frak{F}}_{T}^{\perp}=T_{z_{0}}(T^{\ast}\mathbb{R}^{n})/\mathbb{\frak{F}}_{T}$.
To this linear subspace we associate its dimension
$l_{T}=2d_{T}=\dim (\mathbb{\frak{F}}_{T})$, and also the
following three objects :
\begin{gather}
Q_{T}=p_{2}|\mathbb{\frak{F}}_{T}, \\
\frac{1}{2}q_{T}=p_{2}|\mathbb{\frak{F}}_{T}^{\perp }, \\
C_{Q_{T}}=\{(x,\xi )\in \mathbb{\frak{F}}_{T}\text{ }/\text{
}Q_{T}(x,\xi )=0\}.
\end{gather}
We say that $T\neq 0$ is a period of $d_{x,\xi }\Phi _{t}(z_{0})$ if $%
\mathbb{\frak{F}}_{T}\neq \{0\}$ and that $T$ is a \textbf{total
period} of $d_{x,\xi }\Phi _{t}(z_{0})$ when
$\mathbb{\frak{F}}_{T}=T_{z_{0}}(T^{\ast }\mathbb{R}^{n})$.
\end{definition}
The next condition is inspired by proposition 2.1 of D. Khuat Duy
\cite{KhD} concerning Schr\"{o}dinger operators and will be useful
to separate the contributions of fixed points from those of the
non-trivial periodic trajectories.\medskip
\newline $(H_{5})$\textit{\ For all period }$T$\textit{\ of
}$d_{x,\xi }\Phi_{t}(z_{0})$\textit{\ there exists neighborhoods }$V_{T}$\textit{\ of }$T$%
\textit{\ and }$U_{T}$\textit{\ of }$z_{0}$ \textit{\ such that
$\Phi _{t}(z)\neq z$ for all }$T\in V_{T}$ \textit{and all } $z\in
U_{T}\backslash \{z_{0}\} \cap \Sigma_{E_{c}}$. \medskip

Note that $(H_{5})$ only concerns the dynamics on the energy
surface $\Sigma_{E_{c}}$ and Khuat Duy loc. cit. has shown that
$(H_{5})$ is always satisfied by a Schr\"{o}dinger-type
Hamiltonian $\xi^{2}+V(x)$. We give in the last section two
examples of non-Schr\"{o}dinger operators satisfying $(H_{5})$, to
show that this class is non-empty. We note $\sigma $\textit{\ }the
usual symplectic form on $T_{z_{0}}(T^{\ast } \mathbb{R}^{n})$ and
by $d_{z}^{l}\Phi _{t}$ the derivative of order $l,$ with respect
to initial conditions, evaluated in $(t,z)$. If $u$ is in a vector
space  $V$ we use the notation $u^{l}=(u,...,u)\in V^{l}$.
\begin{definition}
If $k>1$ is the first integer such that $d_{z_{0}}^{k}\Phi _{T}
\neq 0$, we put
\begin{equation}
R_{k}(z)=\frac{1}{k!}\sigma (z,d_{z_{0}}^{k-1}
\Phi_{T}(z^{k-1})),\text{ } z \in \mathbb{\frak{F}}_{T} ,
\end{equation}
and
\begin{equation}
\tilde{R}_{k}=R_{k} | C_{Q_{T}}\cap \mathbb{S}^{2d_{T}-1},
\end{equation}
the restriction of $R_{k}$ to the regular surface $C_{Q_{T}}\cap
\mathbb{S}^{2d_{T}-1}$. Here $\mathbb{S}^{2d_{T}-1}$ is the unit
sphere, where we are working in local coordinates such that
(\ref{formloc}) holds. Finally, let $dL_{Q_{T}}$ be the
Liouville-measure on this surface, i.e. $dL_{Q_{T}}\wedge
dQ_{T}=d\theta$ on $C_{Q_{T}}\cap \mathbb{S}^{2d_{T}-1}$, where
$d\theta$ is the surface measure of $\mathbb{S}^{2d_{T}-1}$.
\end{definition}
The derivatives $d_{z_{0}}^{l}\Phi _{l}$ are computed in section
4.4 and the relation $d_{z_{0}}^{k}\Phi _{T} \neq 0$ can be stated
in terms of resonances, i.e. the arithmetical properties of the
eigenvalues of $p_2$ or, more precisely, of $Q_{T}$.

If $\Theta $ is a cut-off function near the energy $E_{c}$ we have
\begin{equation}
\gamma (E_{c},h)=\mathrm{Tr}(\varphi (\frac{P_{h}-E_{c}}{h})\Theta
(P_{h}))+\mathcal{O}(h^{\infty }),
\end{equation}
as will show proposition \ref{S1(h)=Tr} below. Hence, modulo
$\mathcal{O}(h^{\infty })$, we can write
\begin{gather*}
\gamma (E_{c},h) =\mathrm{Tr}(\psi _{h}^{w}\varphi
(\frac{P_{h}-E_{c}}{h} )\Theta (P_{h}))+\mathrm{Tr}((1-\psi
_{h}^{w})\varphi (\frac{P_{h}-E_{c}}{h}
)\Theta (P_{h})) \\
=\gamma _{1}(z_{0},E_{c},h)+\gamma _{2}(z_{0},E_{c},h),
\end{gather*}
where $\psi $ is a function with compact support near $z_{0}$ such that $%
(H_{5})$ is valid on $\rm{supp}(\psi )$ and the notation $\psi
_{h}^{w}$ stands for the classical Weyl $h$-quantization. Under
the additional hypothesis of having a clean flow, the asymptotics
of $\gamma _{2}(z_{0},E_{c},h)$ is given by the regular trace
formula.

Without any loss of generality we can suppose that
$\rm{supp}(\hat{\varphi})$ is small enough near some non-zero
period $T$ of $d\Phi _{t}(z_{0})$ such that $T$ is the only period
of $d\Phi _{t}(z_{0})$ on $\rm{supp}(\hat{\varphi})$, as follows
by an easy partition of unity argument, the set of periods of the
linearized flow being discrete. For the remaining contribution
$\gamma _{1}(z_{0},E_{c},h)$ we then obtain
\begin{theorem}
\label{Main result Q positive}Under $(H_{1})$ to $(H_{5})$ and the
assumption that $Q_T$ is positive or negative definite, we have
\begin{equation*}
\gamma _{1}(z_{0},E_{c},h)=C(T)\Lambda _{T}(\varphi )+\mathcal{O}(h^{\frac{1}{2}%
}),\text{ as }h\rightarrow 0,
\end{equation*}
where
\begin{equation*}
C(T) =-\frac{1}{2}\frac{\exp (i\frac{\pi
}{4}\mathrm{sgn}(q_{T}))}{|\det (q_{T})|^{\frac{1}{2}}} \exp (i\pi
\frac{d_{T}-1}{2}\mathrm{sign}(Q_{T}))\Gamma (d_{T}),
\end{equation*}
and
\begin{equation*}
\Lambda _{T}(\varphi ) =\frac{1}{(2\pi )^{1+d_{T}}}\left\langle
(t-T-i0)^{-d_{T}},\frac{(t-T)^{d_{T}}}{|\det (\mathrm{Id}-d\Phi
_{t}(z_{0})|^{\frac{1}{2}}}
\hat{\varphi}(t)e^{itp^{1}(z_{0})}\right\rangle .
\end{equation*}
\end{theorem}
In the case of a non definite $Q_T$ we have :
\begin{theorem}
\label{Main result Q non positive}If $n\geq 2$, under $(H_{1})$ to
$(H_{5})$ and if $R_{k}(z)\neq 0$ for all $z\in
C_{Q_{T}}\backslash \{0\}$, we have
\begin{equation*}
\gamma _{1}(z_{0},E_{c},h) = h^{\frac{2d_{T}+k-2}{k}-d_{T}}(K_{T}
\hat{\varphi}(T)\exp(iTp^{1}(z_{0}))+\mathcal{O}(h^{\frac{1}{k}})),\text{
for } h\rightarrow 0,
\end{equation*}
where
\begin{equation*}
K_{T}=\mu _{k}(T)\exp
(i\pi \dfrac{d_{T}-1}{k}\mathrm{sign}(\tilde{R}_{k}))\int\limits_{C_{Q_{T}}%
\mathbb{\cap S}^{2d_{T}-1}}|\tilde{R}_{k}(\theta )|^{-\frac{2d_{T}-2}{k}%
}dL_{Q_{T}}(\theta ),
\end{equation*}
and
\begin{equation*}
\mu _{k}(T)=-\dfrac{1}{k}\Gamma (\dfrac{2d_{T}-2}{k})\frac{
\exp(i\frac{\pi}{4}\mathrm{sgn}(q_{T}))}{|\det(q_{T})|^{\frac{1}{2}}(2\pi)^{d_{T}+1}}.
\end{equation*}
More generally, if $\{ z\in C_{Q_{T}}\backslash \{0\}:
R_{k}(z)=0\}\neq \emptyset$, but if $\nabla Q_{T}(z),\nabla
R_{k}(z)$ are linearly independent on this set, then the same
result holds with
\begin{equation*}
K_{T}=\mu _{k}(T) \exp(i\pi \dfrac{d_{T}-1}{k})
\int\limits_{C_{Q_{T}}\mathbb{\cap
S}^{2d_{T}-1}}(\tilde{R}_{k}(\theta )+i0)^
{-\frac{2d_{T}-2}{k}}dL_{Q_{T}}(\theta ).
\end{equation*}
\end{theorem}
\begin{remark}
\rm{Theorem \ref{Main result Q positive} contains the class of
Schr\"{o}dinger operators $\xi ^{2}+V(x),$ with $V$ non
degenerate, since in this case $p_{2}|\mathrm{Ker}(d\Phi
_{t}(z_{0})-\mathrm{Id})$ is positive definite. Moreover, the
condition $n\geq 2$ in Theorem \ref{Main result Q non positive} is
natural, since, under our assumptions, in dimension $n=1$ we have
$p_2(x_1,\xi_1)=\frac{1}{2}w_1(x_1^2+\xi_1^2)$, or
$\frac{1}{2}w_1(x_1^2-\xi_1^2)$. In the first case $Q_T$ is
definite and in the second case 0 is the only period of the
linearized flow.}
\end{remark}
\section{Oscillatory representation of $\protect\gamma (E_{c},h)$}
We introduce a cut-off function $\Theta \in C_{0}^{\infty
}(]E_{c}-\varepsilon ,E_{c}+\varepsilon \lbrack )$, such that
$\Theta =1$ near $E_{c}$ and $0\leq \Theta \leq 1$ on
$\mathbb{R}$. We then have the decomposition
\begin{eqnarray}
\gamma (E_{c},h) &=&\sum\limits_{\lambda _{j}(h)\in I_{\varepsilon
}}(1-\Theta )(\lambda _{j}(h))\varphi (\frac{\lambda _{j}(h)-E_{c}}{h}%
)+\sum\limits_{\lambda _{j}(h)\in I_{\varepsilon }}\Theta (\lambda
_{j}(h))\varphi (\frac{\lambda _{j}(h)-E_{c}}{h})  \notag \\
&=&\gamma'(E_{c},h)+\gamma''(E_{c},h).
\end{eqnarray}
\begin{proposition}
$\gamma '(E_{c},h)=\mathcal{O}(h^{\infty }),$ when $h\rightarrow
0.\label{S1(h)=Tr}$
\end{proposition}
{\it Proof.} With $\varphi \in \mathcal{S}(\mathbb{R})$ for all
$k\in \mathbb{N}$ there exist $C_{k}$ such that $|x^{k}\varphi
(x)|\leq C_{k}$ on $\mathbb{R}.$ If $N(h)$ is the number of
eigenvalues inside $[E_{c}-\varepsilon,E_{c}+\varepsilon ]\cap
\rm{supp}(1-\Theta )$, then, by Theorem 3.13 of \cite{Rob}, we
have the estimate
\begin{equation*}
|\gamma '(E_{c},h)|\leq N(h)C_{k}|\frac{\lambda
_{j}(h)-E_{c}}{h}|^{-k},\text{ }N(h)=\mathcal{O}(h^{-n}).
\end{equation*}
On the support of $(1-\Theta )$ we have $|\lambda
_{j}(h)-E_{c}|>\varepsilon _{0}>0$, this gives
\begin{equation*}
|\gamma '(E_{c},h)|\leq N(h)C_{N}\varepsilon _{0}^{-k}h^{k}\leq
c_{N}h^{k-n}.
\end{equation*}
Since it is true for all $k\in \mathbb{N}$ the result follows.
\hfill{$\blacksquare$}\medskip

As a consequence, the asymptotics of $\gamma (E_{c},h)$, modulo
$\mathcal{O}(h^{\infty })$, is given by $\gamma''(E_{c},h)$. Now,
by Fourier transform, and inversion, we have that
\begin{equation*}
\Theta (P_{h})\varphi(\frac{P_{h}-E_{c}}{h})
=\frac{1}{2\pi}\int\limits_{\mathbb{R}}\mathrm{Exp}
(\frac{it}{h}P_{h})e^{-i\frac{tE_{c}}{h}}\hat{\varphi}(t)
\Theta(P_{h})dt.
\end{equation*}
The trace of the left-hand side is then exactly $\gamma
''(E_{c},h)$, and we can write
\begin{equation}
\gamma ''(E_{c},h)=\frac{1}{2\pi }\mathrm{Tr}\int\limits_{\mathbb{R}}\Theta (P_{h})\mathrm{%
Exp}(\frac{it}{h}P_{h})e^{-i\frac{tE_{c}}{h}}\hat{\varphi}(t)dt.
\label{Trace S2(h)}
\end{equation}
The formula (\ref{Trace S2(h)}) uses the localized unitary group
$U_{\Theta ,h}(t)=\Theta (P_{h})\mathrm{Exp}(\frac{it}{h}P_{h}).$
A classical result, see e.g. \cite{DUI1}, about this object is
\begin{proposition}
Let $\Lambda$ be the Lagrangian manifold associated to the flow of
$p$
\begin{equation}
\Lambda =\{(t,\tau ,x,\xi ,y,\eta )\in T^{\ast }\mathbb{R}^{2n+1}\text{ }/%
\text{ }\tau =-p(x,\xi ),\text{ }(x,\xi )=\Phi _{t}(y,\eta
)\},\label{lambdaflow}
\end{equation}
then $U_{\Theta ,h}(t)$ is a $h$-Fourier integral operator (or
$h$-FIO) associated to $\Lambda$. More precisely, there exists for
each N a FIO $U_{\Theta ,h}^{(N)}(t)$ with integral kernel in the
H\"{o}rmander class $I(\mathbb{R}\times \mathbb{R}^{n} \times
\mathbb{R}^{n},\Lambda)$ and operators $R_{h}^{(N)}(t)$ on
$L^{2}(\mathbb{R}^{n})$, with uniformly bounded norms for $0<h\leq
1$ and $t$ in a compact subset of $\mathbb{R}$, such that
\begin{equation}
U_{\Theta ,h}(t)=U_{\Theta ,h}^{(N)}(t)+h^{N}R_{h}^{(N)}(t).
\end{equation}
\end{proposition}
We recall the compactly supported cut-off function $\psi=\psi
(x,\xi)$, on whose support $(H_{5} )$ holds. If $\psi_{1}$ is such
that $\psi_{1}\psi =\psi$, with $\rm{supp} (\psi_1)$ small enough,
then by cyclicity of the trace,
\begin{equation*}
\gamma _{2}(z_{0},E_{c},h)=
\frac{1}{2\pi}\mathrm{Tr}\int\limits_{\mathbb{R}}\hat{\varphi}(t)
\psi_{h}^{w}\Theta(P_{h})\mathrm{Exp}(\frac{it}{h}(P_{h}-E_{c}))
\psi_{1,h}^{w}dt+\mathcal{O}(h^{\infty}),
\end{equation*}
and $\psi_{h}^{w}\Theta(P_{h})
\mathrm{Exp}(\frac{it}{h}(P_{h}-E_{c}))\psi_{1,h}^{w}\in
I(\mathbb{R}\times \mathbb{R}^{n} \times \mathbb{R}^{n},\Lambda)$.
After perhaps a local change of variable in $y$, the operator
$\psi_{h}^{w}\Theta(P_{h})
\mathrm{Exp}(\frac{it}{h}P_{h})\psi_{1,h}^{w}$ can be approximated
, modulo an error $\mathcal{O}(h^{N})$, by an $h$-FIO with kernel
\begin{equation*}
\frac{1}{(2\pi h)^{n}}\int\limits_{\mathbb{R}^n} e^{\frac{i}{h}%
(S(t,x,\xi )-\left\langle y,\xi \right\rangle )}b_{N}(t,x,y,\xi
,h)d\xi ,
\end{equation*}
see 4.1 below. Integrating this kernel on the diagonal gives,
modulo terms of order $\mathcal{O}(h^{N})$
\begin{gather}
\gamma _{2}(E_{c},h)=\frac{1}{(2\pi
h)^{n}}\int\limits_{\mathbb{R}}\int\limits_{T^{\ast
}\mathbb{R}^{n}}e^{\frac{i}{h}(S(t,x,\xi )-\left\langle x,\xi
\right\rangle +tE_{c})}a_{N}(t,x,\xi ,h)dtdxd\xi , \label{IO}\\
a_{N}(t,x,\xi ,h)=\hat{\varphi}(t)b_{N}(t,x,x,\xi ,h).
\end{gather}
For a detailed construction we refer to \cite{Rob} or \cite{DUI1}.
\begin{remark}
\rm{Because of the presence of $\Theta (P_{h})$ and
$\hat{\varphi},$ the amplitudes $a_{N}$ are of compact support.
Since we are interested in the main contribution to the trace
formula we note $a$ the amplitude of (\ref{IO}), i.e. $a$ depends
on $(h,N)$.}
\end{remark}
\section{Study of the phase function and of the classical dynamics}
First, we study the nature of the critical points of (\ref{IO}).
After we establish some results on the classical dynamics related
to $p(x,\xi )$ and we compute the Taylor expansion of the phase
function. Resonance-type conditions will naturally occur in the
study of this question.
\subsection{Singularity of the phase}
By Theorem 5.3 of \cite{HOR2}, we can, after perhaps a local
change of variables in $y$, suppose that the flow $\Phi_{t}$, near
$(x_{0},\xi_{0})$ and for $t\in \mathrm{supp}(\hat{\varphi})$
sufficiently small, has a generating function $S(t,x,\eta)$, i.e.
\begin{equation}
(x,\xi )=\Phi _{t}(y,\eta )\Leftrightarrow \left\{
\begin{array}{c}
y=\partial _{\eta }S(t,x,\eta ), \\
\xi =\partial _{x}S(t,x,\eta ),
\end{array}
\right.  \label{Formule generatrice de flot}
\end{equation}
which therefore, by a classical result, satisfies the
Hamilton-Jacobi equation $\partial _{t}S(t,x,\eta )+p(x,\partial
_{x}S(t,x,\eta )=0$. Hence, near $(T,x_{0},\xi_{0})$, the
Lagrangian manifold $\Lambda $ of the flow is parameterized by the
phase function $S(t,x,\eta )-\left\langle y,\eta\right\rangle$.
This choice for the phase is only valid near $(x_{0},\xi_{0})$
when $\xi_{0}\neq 0$, but if $\xi_{0}=0$ we can change the
operator $P_{h}$ by $e^{\frac{i}{h}\left\langle x,\xi_{1}
\right\rangle} P_{h} e^{-\frac{i}{h}\left\langle x,\xi_{1}
\right\rangle}$ with $\xi_{1}\neq 0$. This does not affect the
spectrum and the new operator obtained has symbol
$p(x,\xi-\xi_{1})$ and critical point $(x_0 ,\xi_1)$. With
$(\ref{Formule generatrice de flot})$ a critical point of
$(\ref{IO})$ satisfies
\begin{equation*}
\left\{
\begin{array}{c}
\partial _{x}S(t,x,\xi )=\xi , \\
\partial _{\xi }S(t,x,\xi )=x, \\
\partial _{t}S(t,x,\xi )=-E_{c}.
\end{array}
\right. \Leftrightarrow \left\{
\begin{array}{c}
\Phi _{t}(x,\xi )=(x,\xi {),} \\
p(x,\xi )=E_{c}.
\end{array}
\right.
\end{equation*}
The following lemma is classical and can for example be found in
\cite{KhD1}. Recall that we will often denote points $(x,\xi)$ of
phase space by a single letter $z$.
\begin{lemma}
\label{Lien phase/flot}Let us define $\Psi (t,x,\xi )=S(t,x,\xi
)-\left\langle x,\xi \right\rangle +tE_{c},$ then if $z_{0}$ is
critical point of $\Psi ,$ we have the equivalence $d_{z}^{2}\Psi
(z_{0})\delta z=0\Leftrightarrow d_{z}\Phi _{t}(z_{0})\delta
z=\delta z $, $\forall \delta z\in T_{z_{0}}(T^{\ast
}\mathbb{R}^{n})$, i.e. the degenerate directions of the phase
correspond to fixed points of the linearized flow at $z_0$.
\end{lemma}
If we use Lemma \ref{Lien phase/flot}, we obtain for our phase
function
\begin{corollary}
\label{theo periode flot linéarise}A critical point $(T,x_{0},\xi
_{0})$ of $\Psi (t,x,\xi )$ is degenerate with respect to $(x,\xi
)$ if and only if $T$ is a period of the linearized flow $d_{x,\xi
}\Phi _{t}(x_{0},\xi _{0}).$
\end{corollary}
The next result is also well known from classical mechanics
\begin{lemma} \label{formule flot linéarisé}If $\partial
_{x}p(x_{0},\xi _{0})=\partial _{\xi }p(x_{0},\xi _{0})=0$ then
$d_{x,\xi }\Phi _{t}(x_{0},\xi _{0})$ is the Hamiltonian flow of
the quadratic form $\frac{1}{2}d^{2}p(x_{0},\xi _{0})$ on
$T_{x_{0},\xi _{0}}(T^{\ast }\mathbb{R}^{n})$.
\end{lemma}
Let $T\neq 0$ be a period  of $d\Phi _{t}(z_{0})$. Corollary
\ref{theo periode flot linéarise} shows that we must introduce
\begin{equation}
\Psi (t,x,\xi )=S(t,x,\xi )-\left\langle x,\xi \right\rangle
=(t-T)g(t,x,\xi )+(S(T,x,\xi )-\left\langle x,\xi \right\rangle ).
\label{Decomposition phase en 2 termes}
\end{equation}
This function $g$ is $C^{\infty }$ and satisfies
\begin{equation*}
g(T,x,\xi )=\frac{\partial S}{\partial t}(T,x,\xi )=-p(x,\partial
_{x}S(T,x,\xi )).
\end{equation*}
To simplify the notations we write $g(t,z)=g(t,x,\xi )$ and
\begin{equation}
R(z) =R(x,\xi )=S(T,x,\xi )-\left\langle x,\xi \right\rangle .
\label{Def de R}
\end{equation}
\begin{lemma}
In a neighborhood of $z_{0}$, and near $T,$ the only critical
point, on the energy surface $\Sigma_{E_{c}}$,
of the functions $S(T,x,\xi )-\left\langle x,\xi \right\rangle$ and $%
g(t,x,\xi )$ is $z_{0}$.
\end{lemma}
{\it Proof.} First, $d_{x,\xi }(S(T,x,\xi )-\left\langle x,\xi
\right\rangle )=0$ is equivalent to $\Phi _{T}(x,\xi )=(x,\xi )$.
But near $z_{0}$ and with $(H_{5})$ this can only be satisfied for
$(x,\xi )=z_{0}$. \newline Next, we consider $g(T,x,\xi
)=-p(x,\partial _{x}S(T,x,\xi ))$. here $d_{x,\xi }g(T,x,\xi )=0$
is equivalent to
\begin{equation*}
\left\{
\begin{array}{l}
\partial _{x}p(x,\partial _{x}S(T,x,\xi ))+\partial _{\xi }p(x,\partial
_{x}S(T,x,\xi ))(\partial _{x,x}^{2}S(T,x,\xi ))=0 \\
\partial _{\xi }p(x,\partial _{x}S(T,x,\xi ))(\partial _{x,\xi
}^{2}S(T,x,\xi ))=0.
\end{array}
\right.
\end{equation*}
Clearly $z_{0}$ is critical since $\partial _{x}S(T,x_{0},\xi
_{0})=\xi _{0}$. If $z$ is a critical point of $g(T,x,\xi )$ we
have
\begin{equation}
\left(
\begin{array}{ll}
\mathrm{I}_n & \partial _{x,x}^{2}S(T,x,\xi ) \\
0 & \partial _{x,\xi }^{2}S(T,x,\xi )
\end{array}
\right) \left(
\begin{array}{l}
\partial _{x}p(x,\partial _{x}S(T,x,\xi )) \\
\partial _{\xi }p(x,\partial _{x}S(T,x,\xi ))
\end{array}
\right) =0.  \label{Formule linéaire de Dphi}
\end{equation}
Where $\mathrm{I}_n$ is the identity matrix of order $n$. But
$\Phi _{t}(\partial _{\xi }S(t,x,\xi ),\xi )=(x,\partial
_{x}S(t,x,\xi ))$ and this leads to
\begin{equation*}
d_{x,\xi }\Phi _{t}(\partial _{\xi }S(t,x,\xi ),\xi )\left(
\begin{array}{cc}
\partial _{\xi ,x}^{2}S(t,x,\xi ) & \partial _{\xi ,\xi }^{2}S(t,x,\xi ) \\
0 & \mathrm{I}_n
\end{array}
\right) =\left(
\begin{array}{cc}
\mathrm{I}_n & 0 \\
\partial _{x,x}^{2}S(t,x,\xi ) & \partial _{x,\xi }^{2}S(t,x,\xi )
\end{array}
\right) .
\end{equation*}
Since $d\Phi _{t}$ is an isomorphism Eq.(\ref{Formule linéaire de
Dphi}) imposes
\begin{equation*}
\partial _{x}p(x,\partial _{x}S(T,x,\xi ))=\partial _{\xi }p(x,\partial
_{x}S(T,x,\xi ))=0.
\end{equation*}
In a suitable neighborhood of $z_{0}$ this implies that $(x,\xi
)=z_{0}$, since $z_{0}$ is non-degenerate. With
$(t-T)g(t,z)=S(t,z)-S(T,z),$ for $t\neq T$, the critical points of
$g(t,z)$ are those of $S(t,z)-S(T,z)$, and equation (\ref{Formule
generatrice de flot}) shows that
\begin{gather*}
\Phi _{t}(\partial _{\xi }S(t,x,\xi ),\xi ) =(x,\partial
_{x}S(t,x,\xi))=(x,\partial _{x}S(T,x,\xi ) \\
=\Phi_{T}(\partial_{\xi }S(T,x,\xi ),\xi )=\Phi _{t}
(\partial_{\xi }S(T,x,\xi ),\xi ).
\end{gather*}
We then have by the group law
\begin{equation*}
\Phi _{-t}(\Phi _{T}(\partial _{\xi }S(T,x,\xi ),\xi ))=\Phi
_{T-t}(\partial _{\xi }S(t,x,\xi ),\xi )=(\partial _{\xi
}S(t,x,\xi ),\xi )).
\end{equation*}
For $t\neq T$ and $|t-T|$ small the point $(\partial _{\xi
}S(t,x,\xi ),\xi ))$ would be periodic, with period $(t-T).$ For
$(x,\xi )\in \Sigma_{E_{c}}$, and near $z_{0}$, ($H_{5}$) implies
that $(x,\xi )=z_{0}$. \hfill{$\blacksquare$}
\medskip

The Hessian matrix with respect to $z=(x,\xi )$ of $g$ in
$(T,z_{0})$ satisfies
\begin{equation}
\left\{
\begin{array}{c}
\left\langle \mathrm{Hess}_{z}(g)(T,z_{0})(\delta x,\delta \xi),
(\delta x,\delta \xi )\right\rangle =\left\langle \mathrm{Hess}
(p)(z_{0})B(\delta x,\delta \xi ),B(\delta _{x},\delta \xi
)\right\rangle ,\\
B=\left(
\begin{array}{cc}
\mathrm{I}_n & 0 \\
\partial _{x,x}^{2}S(T,z_{0}) & \partial _{x,\xi }^{2}S(T,z_{0})
\end{array}
\right) ,
\end{array}
\right.
\end{equation}
with $B$ non singular, as seen before. This proves the relation
\begin{equation}
\mathrm{Hess}_{z}(g)(T,z_{0})=\text{
}^{t}B\mathrm{Hess}(p)(z_{0})B.
\end{equation}
In the case of a total period $T$ of $d\Phi _{t}(z_{0})$ we obtain
\begin{proposition}
\label{Hessiennefullperiod}If $d\Phi _{t}(z_{0})$ is totally
periodic, with period $T$, then the function $g(t,x,\xi )$
satisfies $\mathrm{Hess}_{z}(g)(T,z_{0})=\mathrm{Hess}(p)(z_{0}).$
\end{proposition}
\subsection{The linearized flow in $z_{0}$}
Up to a permutation of coordinates we can assume that
\begin{equation}
p_{2}(x,\xi )=\frac{1}{2}(\sum_{j=1}^{k}w_{j}(x_{j}^{2}+\xi
_{j}^{2})+\sum_{j=k+1}^{n}w_{j}(x_{j}^{2}-\xi _{j}^{2})).
\label{truc}
\end{equation}
The flow of $p_{2}$, viewed as an element of $\mathrm{End}
(T_{0}(T^{\ast }\mathbb{R}^{n}))\simeq
\mathrm{End}(\mathbb{R}^{2n})$, is
\begin{equation*}
\mathrm{Exp}(tH_{p_{2}})(x,\xi )=A(t)\left(
\begin{array}{l}
x \\
\xi
\end{array}
\right) ,\text{ }A(t)=\left(
\begin{array}{cccc}
a(t) & 0 & c(t) & 0 \\
0 & b(t) & 0 & d(t) \\
-c(t) & 0 & a(t) & 0 \\
0 & d(t) & 0 & b(t)
\end{array}
\right) ,
\end{equation*}
where $(x,\xi )=(x',x'',\xi',\xi'')$, $x',\xi' \in \mathbb{R}^k$,
$x'',\xi'' \in  \mathbb{R}^{n-k}$, and
\begin{equation*}
\left\{
\begin{array}{l}
a(t)=\mathrm{diag}(\cos (w_{i}t)),\text{ }
b(t)=\mathrm{diag}(\mathrm{ch}(w_{i}t)), \\
c(t)=\mathrm{diag}(\sin (w_{i}t)), \text{ }
d(t)=\mathrm{diag}(\mathrm{sh}(w_{i}t)),
\end{array}
\right.
\end{equation*}
and "diag" means diagonal matrix. In the following we work on the
subspace $\{x''=\xi''=0\}$ obtained by projection on the periodic
variables. Let be $I$ a subset of $\{1,...,n\}$ with $l$ elements,
$l>1.$ The existence of a non trivial closed trajectory of
dimension $l$ imposes that there exists a $c\in \mathbb{R}^{\ast}$
such that
\begin{equation}
\forall i\in I,\text{ }\exists n_{i}\in (%
\mathbb{Z}^{\ast })^{l}\text{ }:\text{ }w_{i}=cn_{i}.
\label{Dépendance sur Z des fréquences}
\end{equation}
\begin{remark}
\rm{Let $M(w)=\{k\in \mathbb{Z}^{n}\text{ }/\text{ }\left\langle
k,w\right\rangle =0\}$ be the $\mathbb{Z}$-module of resonances of
the vector $w=(w_{1},...,w_{n})$. The relations (\ref{Dépendance
sur Z des fréquences}), for $l>1$ lead to resonances, since
$n_{i}w_{i}-n_{j}w_{j}=0$, but a resonant system can have no
periodic trajectories of dimension greater than 1, as is shown by
$(\sqrt{2}+\sqrt{3},\sqrt{2}-\sqrt{3},\sqrt{2}).$}
\end{remark}
Since we are interested by the periods of $d\Phi _{t}(z_{0})$ we
define (assuming (\ref{truc}))
\begin{equation*}
Q_{per}(x,\xi )
=\frac{1}{2}\sum\limits_{j=1}^{k}w_{j}(x_{j}^{2}+\xi _{j}^{2}),
\end{equation*}
and also
\begin{equation*}
Q^{+}(x,\xi
)=\sum\limits_{j=1}^{k_{1}}\frac{w_{j}}{2}(x_{j}^{2}+\xi
_{j}^{2}),\text{ }w_{i}>0,\end{equation*}
\begin{equation*}
Q^{-}(x,\xi)=\sum\limits_{j=k_{1}+1}^{k}\frac{|w_{j}|}{2}(x_{j}^{2}+
\xi _{j}^{2}),\text{ }w_{i}<0,
\end{equation*}
so that $Q_{per}(x,\xi )=Q^{+}(x,\xi )-Q^{-}(x,\xi )$. Let $P_{1}$
and $P_{2}$ be the linear subspaces obtained by projecting
orthogonally on the effective variables of $Q^{+}$, $Q^{-}$.
\begin{proposition}
\label{Hess g, restric definie}If we have $\mathbb{\frak{F}}%
_{T}\subset P_{1}$ or $\mathbb{\frak{F}}_{T}\subset P_{2},$ then $\mathrm{%
Hess}(g)(T,z_{0})_{|\mathbb{\frak{F}}_{T}}$ is respectively
positive or negative definite.
\end{proposition}
{\it Proof.} Choosing coordinates as in (\ref{truc}), we see that
for $v\in \frak{F}_T$
\begin{equation*}
A(T) \left(
\begin{array}{cc}
\partial _{x,\xi }^{2}S(T,z_{0}) & \partial _{\xi ,\xi }^{2}S(T,z_{0}) \\
0 & \mathrm{I}_{n}
\end{array}
\right) v =\left(
\begin{array}{cc}
\mathrm{I}_{n} & 0 \\
\partial _{x,x}^{2}S(T,z_{0}) & \partial _{x,\xi }^{2}S(T,z_{0})
\end{array}
\right)v .
\end{equation*}
With $\textrm{dim}(\mathbb{\frak{F}}_{T})=2d_{T}$ we can choose
our coordinates such that
\begin{equation*}
\left(
\begin{array}{cc}
a(T) & 0 \\
0 & b(T)
\end{array}
\right) \partial _{x,\xi }^{2}S(T,z_{0})=\left(
\begin{array}{cc}
\left(
\begin{array}{ll}
\mathrm{I}_{d_{T}} & 0 \\
0 & \ast
\end{array}
\right) & 0 \\
0 & b(T)
\end{array}
\right) \partial _{x,\xi }^{2}S(T,z_{0})=\mathrm{I}_{n}.
\end{equation*}
Elementary considerations show that
\begin{equation}
\partial _{x,\xi }^{2}S(T,z_{0})=\left(
\begin{array}{cc}
\mathrm{I}_{d_{T}} & 0 \\
0 & \ast
\end{array}
\right) ,\text{ }\partial _{x,x}^{2}S(T,z_{0})=\left(
\begin{array}{cc}
0_{d_{T}} & 0 \\
0 & \ast
\end{array}
\right) ,
\end{equation}
where $\ast$ designs matrix blocs which are irrelevant for the
present discussion. Hence, by restriction to
$\mathbb{\frak{F}}_{T}$ : $
(\mathrm{Hess}(g)(T,z_{0}))_{|\mathbb{\frak{F}}_{T}}
=(\mathrm{Hess}(p)(z_{0}))_{|\mathbb{\frak{F}}_{T}}$.\hfill{$\blacksquare$}
\subsection{Taylor series of the flow near $z_{0}$}
We start by the general case of an autonomous system near an
equilibrium. Let be $\Phi _{t}$ the flow of a $C^{\infty }$ vector
field $X$ on $\mathbb{R}^{n}$ with coordinates
$z=(z_{1},...,z_{n})$ and let $z_{0}$ a fixed point of $\Phi
_{t}$. We denote by $A(z_{0})$ the matrix of the linearization of
$X$ in $z_{0}$
\begin{equation}
A(z_{0})=(\frac{\partial X^{i}}{\partial z_{k}})_{(i,k)}(z_{0}),
\end{equation}
and we recall that $d\Phi_{t}(z_{0})=\mathrm{Exp}(tA(z_{0}))$.
Here, and in the following, the derivatives $d$ will be taken with
respect to $z,$ we denote by $d^{k}f$ the $k$-th derivative of $f$
regarded as a multi-linear form on the $k$-fold product
$\mathbb{R}^{n}\times ...\times \mathbb{R}^{n}$, and
$d^{k}f(z_{0})$ or $d_{z_{0}}^{k}f$ this derivative evaluated in
$z_{0}$. As an example, we compute the second derivative of the
flow in $z_{0}$ :
\begin{gather}
\partial _{z_{i},z_{j}}^{2}(\frac{d}{dt}\Phi _{t}(z))=\sum_{k,l=1}^{n}\frac{%
\partial ^{2}X}{\partial z_{k}\partial z_{l}}(\Phi _{t}(z))\frac{\partial
\Phi _{t}^{l}(z)}{\partial z_{j}}\frac{\partial \Phi
_{t}^{k}(z)}{\partial z_{i}}\notag\\%
+\sum_{k=1}^{n}\frac{\partial X}{\partial z_{k}}(\Phi _{t}(z))\frac{%
\partial ^{2}\Phi _{t}^{k}(z)}{\partial z_{i}\partial z_{j}}. \label{equation D2phi}
\end{gather}
Hence at the point $z_{0}$ we obtain
\begin{equation*}
\frac{d}{dt}(\partial _{z_{i},z_{j}}^{2}\Phi _{t}(z_{0}))=\sum_{k,l=1}^{n}%
\frac{\partial ^{2}X}{\partial z_{k}\partial
z_{l}}(z_{0})\frac{\partial
\Phi _{t}^{l}(z_{0})}{\partial z_{j}}\frac{\partial \Phi _{t}^{k}(z_{0})}{%
\partial z_{i}}+A(z_{0})\partial _{z_{i},z_{j}}^{2}\Phi _{t}(z_{0}).
\end{equation*}
Let us write $\mathrm{Hess}(X)(z_{0})$ for the vector valued
Hessian of $X$ evaluated in $z_{0}$. Interpreting (\ref{equation
D2phi}) as an inhomogeneous system of equations for
$\partial_{z_{i},z_{j}}\Phi _{t}(z_{0})$ we obtain, since
$\Phi_{0} =\mathrm{Id}$ and therefore $d_{z,z}^{2}\Phi
_{0}(z_{0})\equiv 0$, that
\begin{equation}
d^{2}\Phi _{t}(z_{0})(z,z)=d\Phi
_{t}(z_{0})\int\limits_{0}^{t}d\Phi
_{-s}(z_{0})\mathrm{Hess}(X)(z_{0})(d\Phi _{s}(z_{0})(z),d\Phi
_{s}(z_{0})(z))ds.
\end{equation}
Now, let us assume that $z_{0}$ is the origin, we generalize as
follows :
\begin{proposition} \label{prop deg k}
If $X(z)=Az+X_{k}(z)+\mathcal{O}(||z||^{k+1}),$ where $X_{k}(z)$
is homogeneous of degree $k>2$ in $z$, then the flow $\Phi _{t}$
of $X$ satisfies
\begin{gather*}
d^{j}\Phi _{t}(0)=0,\text{ }\forall j\in \{2,...,k-1\}, \\
d^{k}\Phi _{t}(0)(z,...,z)=d\Phi _{t}(0)\int\limits_{0}^{t}d\Phi
_{-s}(0)d^{k}X(0)(d\Phi _{s}(0)(z),...,d\Phi _{s}(0)(z))ds.
\end{gather*}
\end{proposition}
{\it Proof.} The first result is trivial. At the order $k$ we have
for $|\alpha | = k$
\begin{equation*}
\frac{d}{dt}((\frac{\partial }{\partial z})^{\alpha }\Phi _{t}(z))=(\frac{%
\partial }{\partial z})^{\alpha }(X(\Phi _{t}(z)).
\end{equation*}
Hence, for $z=0$ we simply have
\begin{equation}
\frac{d}{dt}((\frac{\partial }{\partial z})^{\alpha }\Phi
_{t}(z))_{|z=0}=d^{\alpha }X(0)(d\Phi _{t}(z_{0}),...,d\Phi
_{t}(z_{0}))+A((\frac{\partial }{\partial z})^{\alpha }\Phi
_{t}(z))_{|z=0}. \label{Formule de p-flot}
\end{equation}
Eq. (\ref{Formule de p-flot}) is linear of the form
$d_{t}u_{\alpha }(t)=f_{\alpha }(t)+A u_{\alpha }(t)$, and by
integration, with the initial condition $d^{k}\Phi _{0}(0)=0,$ we
obtain
\begin{equation*}
((\frac{\partial }{\partial z})^{\alpha }\Phi _{t})(0)=d\Phi
_{t}(0)\int\limits_{0}^{t}d\Phi _{-s}(0)d^{\alpha }X(0)(d\Phi
_{s}(0),...,d\Phi _{s}(0))ds,
\end{equation*}
and the result holds by linearity. \hfill{$\blacksquare$}\medskip

A more general result is
\begin{theorem}
\label{TheoFormule de récurence du flot}Let be $z_{0}$ an
equilibrium point of $X$ and $\Phi _{t}$ the flow of $X.$ For all
$m\in \mathbb{N}^{\ast }$, there exists a polynomial map $P_{m}$,
of degree at most $m$, such that
\begin{equation}
d^{m}\Phi _{t}(z_{0})(z^{m})=d\Phi
_{t}(z_{0})\int\limits_{0}^{t}d\Phi _{-s}(z_{0})P_{m}(d\Phi
_{s}(z_{0})(z),...,d^{m-1}\Phi _{s}(z_{0})(z^{m-1}))ds.
\label{Formule de récurence du flot}
\end{equation}
In addition $P_{m}$ is uniquely determined by the $m$-jet of $X$
in $z_{0}$.
\end{theorem}
{\it Proof.} For $m=1$, $d\Phi_t(z_{0})$ is determined by the
operator $A(z_{0})$, i.e. by the $1$-jet of $X$. We note $x^l\in
(\mathbb{R}^n)^l$ the image of $x$ under the diagonal mapping,
with the same convention for any vector. If $f$ and $g$ are smooth
we obtain
\begin{equation*}
d^{m}(fog)(x_{0})(x^{m})=\sum\limits_{j=1}^{m}\sum\limits_{\alpha
\in I_{j}}c_{\alpha }d^{j}f(g(x_{0}))((dg(x_{0})(x))^{\alpha
_{j_{1}}},...,(d^{m}g(x_{0})(x^{m}))^{\alpha _{j_{m}}}),
\end{equation*}
with $c_{\alpha }\in \mathbb{N}$, $I_{j}=\{\alpha \in
\mathbb{N}^{m}$ $/$ $\sum\limits_{k=1}^{m}k\alpha _{j_{k}}=j\}$.
Since $z_{0}$ is a fixed point, we have
\begin{equation*}
d^{m}(Xo\Phi
_{t})(z_{0})(z^{m})=\sum\limits_{j=1}^{m}\sum\limits_{\alpha \in
I_{j}}c_{\alpha }d^{j}X(z_{0})((d\Phi _{s}(z_{0})(z))^{\alpha
_{j_{1}}},...,(d^{m}\Phi _{s}(z_{0})(z^{m}))^{\alpha _{j_{m}}}).
\end{equation*}
For $Y=(Y_{1},...,Y_{m})$, we can define
\begin{equation}
\label{def P_m} P_{m}(Y)=\sum\limits_{j=1}^{m}\sum\limits_{\alpha
\in I_{j}}c_{\alpha }d^{j}X(z_{0})(Y_{1}^{\alpha
_{j_{1}}},....,Y_{m}^{\alpha _{m}})-dX(z_{0})(Y_{m}),
\end{equation}
this leads to the differential equation, operator valued
\begin{equation*}
\frac{d}{dt}(d^{m}\Phi _{t}(z_{0}))(z^{m})=A(z_{0})d^{m}\Phi
_{t}(z_{0})(z^{m})+P_{m}(d\Phi _{s}(z_{0})(z),...,d^{m-1}\Phi
_{s}(z_{0})(z^{m-1})).
\end{equation*}
With the initial condition $d^{m}\Phi _{0}(z_{0})=0,$ we obtain
that the solution is given by (\ref{Formule de récurence du
flot}). Moreover, Eq.(\ref{def P_m}) shows that $P_{m}$ is
completely determined by the derivatives of order less or equal
than $m$ of $X.$ \hfill{$\blacksquare$}
\subsection{Application to Hamiltonian systems} Proposition \ref{prop deg k} applied to
the flow of $H_{p}$ shows that
\begin{equation*}
d^{2}\Phi _{t}(z_{0})(z,z)=d\Phi
_{t}(z_{0})\int\limits_{0}^{t}d\Phi
_{-s}(z_{0})\mathrm{Hess}(H_{p})(z_{0})(d\Phi _{s}(z_{0})(z),d\Phi
_{s}(z_{0})(z))ds.
\end{equation*}
We now consider more closely $d^{2}\Phi _{t}(z_{0})$ for $t=T$, a period of $%
d\Phi _{t}(z_{0})$. We introduce the following terminology
\begin{definition}
$w\in (\mathbb{R}^{\ast })^{n}$ is pseudo-resonant to the order
$l\in \mathbb{N}$ if
\begin{equation}
\exists \text{ }(i_{1},...,i_{l}) \in \{1,...,n\},(\varepsilon
_{i_{1}},...,\varepsilon _{i_{l}})\in \{-1,1\}\text{ such that
}\sum\limits_{j=1}^{l}\varepsilon _{i_{j}}w_{i_{j}}=0.
\label{pseudoresonances}
\end{equation}
\end{definition}
\begin{remark}
\rm{This notion is weaker than the usual resonance condition since
for $l$ even there always exists a pseudo-resonance. For example,
for $l=4$ we have $(w_{i} -w_{i})\pm  (w_{j} -w_{j})=0 $, although
$w$ can be non-resonant at the order 4. We also observe that all
resonances of order 3 are pseudo-resonances.}
\end{remark}
In term of the frequencies $w_{i}$, we then have
\begin{theorem}
\label{theo non auto D2}If the frequencies $w$ satisfy no
pseudo-resonance relation of order 3 and if $T$ is a total period
of $d\Phi _{t}(z_{0})$, we have $d^{2}\Phi _{T}(z_{0})=0.$
\end{theorem}
{\it Proof.} Under the condition $(H_{3})$, $d^{2}\Phi
_{t}(z_{0})$ can be expressed as a linear combination of integrals
of the elementary functions $s\mapsto \exp (\pm
is(w_{i}+\varepsilon _{1}w_{j}+\varepsilon _{2}w_{k}))$ with
$\varepsilon _{j}=\pm 1$. Hence, to determine $d^{2}\Phi
_{T}(z_{0})$ we must compute
\begin{equation*}
\int\limits_{0}^{T}\exp (\pm is(w_{i}+\varepsilon
_{1}w_{j}+\varepsilon _{2}w_{k}))ds,
\end{equation*}
but under the assumptions of the theorem all these integrals are
$0$. \hfill{$\blacksquare$}\medskip

For all $l\in \mathbb{N}^{\star}$ let be $M_{l}(W) =
\{k\in\mathbb{Z}^{2n}\text{ / }\langle k,(w,w)\rangle =0, \text{
}|k|=l\}$ the $\mathbb{Z}$-module of resonances of order $l$. In
the presence of resonances we can say the following :
\begin{proposition} For $k\in M_{3}(W)$ let be
$\tilde{k}=(|k_{1}|,...,|k_{2n}|).$ If $\forall k\in M_{3}(W)$ we
have $\dfrac{\partial ^{3}p}{\partial z^{\tilde{k}}}(z_{0})=0$ and
if $T$ is a total period of $d_{z_{0}}\Phi _{t}$ then $d^{2}\Phi
_{T}(z_{0})=0$.
\end{proposition}
The proof is trivial when going back to the proof of Theorem
\ref{theo non auto D2}. \hfill{$\blacksquare $}\medskip \newline
Let  $f^{\ast }$ be the pullback by a map $f$. Then by proposition
\ref{prop deg k} we have :
\begin{corollary}
\label{relevant} For a Hamiltonian system with the equilibrium
point $z_{0}$ and such that $d_{z_{0}}^{j}p=0,$ $\forall j\in
\{3,...,k-1\}$, we obtain
\begin{gather*}
d^{j}\Phi _{t}(z_{0})=0,\text{ }\forall t,\text{ }\forall j\in
\{2,...,k-2\},\\
d_{z_{0}}^{k-1}\Phi _{t}(z^{k-1})=d\Phi
_{t}(z_{0})\int\limits_{0}^{t}d\Phi _{-s}(z_{0})\left( d\Phi
_{s}(z_{0})^{\ast }(d_{z_{0}}^{k-1}H_{p})\right)(z^{k-1}) ds.
\end{gather*}
\end{corollary}
And Theorem \ref{theo non auto D2} generalizes trivially to the
order $k$ under the conditions of Corollary \ref{relevant}. More
precisely, if $k$ is odd and if there is no pseudo-resonance of
order $k$ then, under the assumptions of Corollary \ref{relevant},
we have $d_{z_{0}}^{k-1}\Phi _{t}=0$.
\subsection{Relation between the phase and the flow}
Like in the preceding section we consider a Hamiltonian function
$p$ with total period $T$ for the linearized flow, satisfying,
near 0
\begin{equation}
p(z)=E_{c}+p_{2}(z)+\mathcal{O}(||z||^{k}). \label{Conditions
ordre p}
\end{equation}
We recall that $z=(x,\xi)$ and
$z^{k}=(z,...,z)\in\mathbb{R}^{2nk}$. By Taylor, we have
\begin{equation*}
\Phi _{T}(z)=z+\frac{1}{(k-1)!}d^{k-1}\Phi
_{T}(0)(z^{k-1})+\mathcal{O}(||z||^{k}).
\end{equation*}
Under these conditions we can write the generating function at
time $T$ as
\begin{equation}
S(T,x,\xi )=\left\langle x,\xi \right\rangle +R_{k}(x,\xi
)+R_{k+1}(x,\xi ),
\end{equation}
where $R_{k}$ is homogeneous of degree $k$ and $R_{k+1}$ is the
remainder of the Taylor expansion. Let $J$ be the matrix of the
standard symplectic form on $T^{\ast}\mathbb{R}^{n}$. The relation
between the phase function and the flow of $p$ is given by
\begin{proposition}
\label{TheopourRk}Under conditions (\ref{Conditions ordre p}) the
$(k$-$1)$-st derivative of the flow $\Phi _{t}$, at time $t=T$,
equals :
\begin{equation}
d^{k-1}\Phi _{T}(0)((x,\xi )^{k-1}) = -(k-1)!J\nabla R_{k}(x,\xi)
\end{equation}
In addition we have $R_{j}(x,\xi )=0,$ $3\leq j<k,$ and
\begin{equation}
R_{k}(x,\xi )=\frac{1}{k!}\int\limits_{0}^{T}\sigma ((x,\xi
),d\Phi _{-s}(0)\left( d\Phi _{s}(0)^{\ast
}(d_{0}^{k-1}H_{p})\right) (x,\xi )^{k-1})ds.\label{formule Rk}
\end{equation}
\end{proposition}
{\it Proof.} With equation (\ref{Formule generatrice de flot}) we
obtain
\begin{equation*}
\Phi _{T}(x+\partial _{\xi }R_{k}+\partial _{\xi }R_{k+1},\xi
)=(x+\partial _{\xi }R_{k},\xi )+\frac{1}{(k-1)!}d^{k-1}\Phi
_{T}(0)((x,\xi )^{k-1})+\mathcal{O}(||(x,\xi )||^{k}),
\end{equation*}
by identification of homogeneous terms we have successively
\begin{gather*}
(-\partial _{\xi }R_{k},\partial _{x}R_{k})(x,\xi )=-J\nabla R_{k}(x,\xi )=%
\frac{1}{(k-1)!}d^{k-1}\Phi _{T}(0)((x,\xi )^{k-1}),\\
R_{k}(x,\xi )=\frac{1}{k!}\left\langle (x,\xi ),-J\nabla
R_{k}(x,\xi )\right\rangle =\frac{1}{k!}\sigma ((x,\xi
),d^{k-1}\Phi _{T}(0)((x,\xi )^{k-1})),
\end{gather*}
where the last result holds by homogeneity. Corollary
\ref{relevant} then implies (\ref{formule Rk}).
\hfill{$\blacksquare$}
\section{Normal forms of the phase function}
In this section we derive suitable normal forms for $\Psi (t, x,
\xi) = S(t, x, \xi ) - \langle x, \xi \rangle $ in our oscillatory
integral representation of $\gamma _{2}(E_c,h)$. We recall the
decomposition
\begin{equation*}
\Psi (t, x , \xi ) = (t - T ) g(t, x , \xi ) + R(x,\xi ) ,
\end{equation*}
cf. formulas (\ref{Decomposition phase en 2 termes}) and (\ref{Def
de R}). In the micro-local neighborhood of $z_0 = (x_0 , \xi _0 )
$ we are interested in, the only critical point of $R $ and of
$g(t, \cdot ) $ is $z_0 $ for $t $ close to $T $ and, moreover,
$z_0 $ is non-degenerate for the latter.

A further very important simplifying assumption we will make for
the moment is that, until further notice, {\it $T $ is a total
period of $d\Phi _t (z_0 )$}. We will show in section 6.3 below
how to relax this assumption. If $T $ is such a total period, then
clearly $R(z) = \mathcal{O}(||z||^3 ) $. This can be made more
precise
\begin{lemma} \label{lemmephase}If near $z_{0}$ the
function $p$ satisfies $(H_{3})$ and  condition (\ref{Conditions
ordre p}) then for $t$ near $T$ there exist a non degenerate
quadratic form $Q_{t}(x,\xi )$ such that $Q_{T}(x,\xi
)=p_{2}(x,\xi)$ and
\begin{equation*}
S(t,x,\xi )-\left\langle x,\xi \right\rangle =(t-T)\left(
Q_{t}(x,\xi )+h(t,x,\xi )\right) +R(x,\xi ),
\end{equation*}
with $R(x,\xi )=\mathcal{O}(||(x,\xi )||^{k})$ and $h(t,x,\xi
)=\mathcal{O}(||(x,\xi )||^{k})$ uniformly in $t$.
\end{lemma}
{\it Proof.} On replacing $t-T$ by $t$ we can write $\Psi = R(z) +
t G(t,z) $ with $G(t,z)=g(t+T,z)$. By a second order Taylor
expansion around $z_0 $ and proposition \ref{Hessiennefullperiod}
we have that $G(t, z ) = Q_t (z) + h(t,z)$, with
$Q_{0}(z)=p_{2}(z)$ and $h (t,z) = \mathcal{O}(||z||^3 )$. Now by
proposition \ref{TheopourRk}, $R(z) = \mathcal{O}(||z||^k)$, given
that $p$ satisfies (\ref{Conditions ordre p}), and since $\Phi _T
(z )=z + \mathcal{O}(||z||^{k - 1 } )$, we have that $S(T, x, \xi
)=\langle x,\xi \rangle + \mathcal{O}(||(x, \xi )||^k ) $.
Therefore $h(t,z)=\mathcal{O}(||z||^{k})$ uniformly in $t$ and the
lemma follows. \hfill{$\blacksquare$}
\medskip \newline
\textbf{Reduction of the phase with respect to
}$C_{Q_{T}}$\textbf{.} \medskip

Let us suppose $S(T,x,\xi )$ contains effectively some terms of
order $k$. We write, as before, $R(z)=R_{k}(z)+R_{k+1}(z)$, where
$R_{k}$ is the homogeneous component of degree $k$ of $R$ and
$R_{k+1}$ is the remainder of the Taylor series. The following
lemma is useful for any perturbation by a function of odd degree.
\begin{lemma}
\label{funnylemma}Let be $Q$ a non degenerate quadratic form on
$\mathbb{R}^{n},$ $n>3,$ with inertia indices greater than 2. For
all odd continuous function $R$, $Q$ and $R$ have a common zero on
$\mathbb{S}^{n-1}$.
\end{lemma}
{\it Proof.} Up to a linear change of coordinates we can assume
that $Q(x) =||x_{1}||^{2}- ||x_{2}||^{2}$, with $x_{1} \in
\mathbb{R}^{p},\text{ }x_{2}\in \mathbb{R}^{q}, \text{ }p,q\geq
2,\text{ }p+q=n.$ Now the cone of the zeros of $Q$ is invariant
under isometries of the subspaces $(x_{1},0)$ and $(0,x_{2}).$ By
rotating around the origins there
exist a continuous curve $\gamma _{1}$ inside the cone of $Q$ mapping $%
(x_{1},x_{2})$ to $(-x_{1},x_{2})$ and a curve $\gamma _{2}$ mapping $%
(-x_{1},x_{2})$ to $(-x_{1},-x_{2})$, inside the cone. If $\gamma
=\gamma _{1}.\gamma _{2}$ is the union of the two previous curves,
the function $R(\gamma )$ gives the result by continuity since $R$
is odd. \hfill{$\blacksquare$}
\begin{remark}
\rm{A consequence of Lemma \ref{funnylemma} is that the set
$C_{Q_T} \cap C_{R_k} \cap \mathbb{S}^{n-1}$ is not empty when the
function $R_{k}$ is non-zero and odd and $Q_T$ is non-definite.}
\end{remark}
We choose polar coordinates $z=(x,\xi )=r\theta $, $\theta \in \mathbb{S%
}^{2n-1}(\mathbb{R}).$ These coordinates will perform a ''blow-up'' of $%
\mathbb{R}\times (T^{\ast }\mathbb{R}^{n}\backslash \{0\})$. In
general one uses the projective space
$\mathbb{P}_{2n-1}(\mathbb{R})$, but here, since the singularities
are carried by the conic set of the zeros of $Q_{T}$, it is
convenient to use the sphere $\mathbb{S}^{2n-1}$. For any function
$f$, positively homogeneous on $\mathbb{R}^{n}$, we note $
C_{f}=\{x\in \mathbb{R}^{n}\text{ }/\text{ }f(x)=0\}$, the conic
set of the zeros of $f$. Finally, $g\simeq h$ means that
applications $g$ and $h$ are conjugated by a local diffeomorphism.
\begin{lemma}
\label{Lemme Carte1}If $\theta _{0}\in \mathbb{S}^{2n-1}$ with
$\theta
_{0}\notin C_{Q_{0}}$, there exists a system of local coordinates $\chi $%
, near $(t,r,\theta )=(0,0,\theta _{0}),$ such that $\Psi \simeq
\chi _{0}\chi _{1}^{2}$ in a neighborhood of $(\chi _{0},\chi
_{1})=(0,0).$
\end{lemma}
{\it Proof.} Using the notations of Lemma \ref{lemmephase} we
have, in polar coordinates,
\begin{equation*}
\Psi (t,z) \simeq r^{2}(tQ_{t}(\theta )+tr^{k-2}h_{1}(t,r,\theta
)+r^{k-2}R_{k}(\theta )+r^{k-1}\tilde{R}(r,\theta )),
\end{equation*}
with $\tilde{R} \in C^{\infty }(\mathbb{R}_{+}\times
\mathbb{S}^{2n-1})$ and $h_{1} \in  C^{\infty }(\mathbb{R\times
R}_{+}\times \mathbb{S}^{2n-1})$. We choose new coordinates
\begin{gather*}
\chi _{0}(t,r,\theta ) =tQ_{t}(\theta )+r^{k-2}(th_{1}(t,r,\theta
)+R_{k}(\theta )+r\tilde{R}(r,\theta )), \\
(\chi _{1},\chi _{2},...,\chi _{2n})(t,r,\theta ) =(r,\theta ).
\end{gather*}
Hence $|\frac{D\chi }{D(t,r,\theta )}|(0,0,\theta
_{0})=|Q_{0}(\theta _{0})|\neq 0$ and in the new system of
coordinates we have $\Psi (t,z)=(\chi _{0}\chi
_{1}^{2})(t,r,\theta )$ as required. \hfill{$\blacksquare$}
\begin{lemma}
\label{Lemme Carte2}If $\theta _{0}\in C_{Q_{0}}$ with $\theta
_{0}\notin C_{R_{k}}$ there exists a system of local coordinates $\chi $%
, near $(t,r,\theta )=(0,0,\theta _{0}),$ such that $\Psi \simeq
\chi _{0} \chi _{2}\chi _{1}^{2}\pm \chi _{1}^{k}$ in a
neighborhood of $(\chi _{0},\chi _{1},\chi _{2})=(0,0,0).$
\end{lemma}
{\it Proof.} As in Lemma \ref{Lemme Carte1} we write
\begin{equation*}
\Psi (t,z)=tr^{2}(Q_{t}(\theta )+r^{k-2}h_{1}(t,r,\theta
))+r^{k}(R_{k}(\theta )+r\tilde{R}(r,\theta )).
\end{equation*}
Since $R_{k}(\theta _{0})\neq 0$ in a suitable neighborhood of
$(\theta _{0},0)$ we have $R_{k}(\theta)+r\tilde{R}(r,\theta )\neq
0.$ Up to a permutation of the variables $\theta _{i},$ we can
suppose that $\frac{\partial Q_{0}}{\partial \theta _{1}}(\theta
_{0})\neq 0.$ In a neighborhood of $(t,r,\theta )=(0,0,\theta
_{0}),$ we choose coordinates
\begin{equation*}
\begin{array}{l}
\chi _{0}(t,r,\theta )=t|R_{k}(\theta )+r\tilde{R}(r,\theta )|^{-\frac{2}{k}%
}, \\
\chi _{1}(t,r,\theta )=r|R_{k}(\theta )+r\tilde{R}(r,\theta
)|^{\frac{1}{k}},\\
\chi _{2}(t,r,\theta )=Q_{t}(\theta )+r^{k-2}h_{1}(t,r,\theta ), \\
(\chi _{3},...,\chi _{2n})(t,r,\theta )=(\theta _{2},...,\theta
_{2n-1}).
\end{array}
\end{equation*}
The corresponding Jacobian, at the point $(0,0,\theta _{0})$, is
\begin{equation}
|\frac{D\chi }{D(t,r,\theta )}(0,0,\theta _{0})|=|R_{k}(\theta _{0})|^{-\frac{1}{k}}|%
\frac{\partial Q_{0}}{\partial \theta _{1}}(\theta _{0})|\neq 0.
\end{equation}
In these new coordinates the phase is $\Psi (t,r,\theta )=(\chi
_{0}\chi _{2}\chi _{1}^{2}+\chi _{1}^{k})(t,r,\theta ).$
\hfill{$\blacksquare$}
\begin{lemma}
\label{Lemme Carte3}If $\theta _{0}\in C_{Q_{0}}\cap C_{R_{k}}$
and if $\nabla Q_{0}(\theta _{0}),\nabla R_{k}(\theta _{0})$ are
linearly independent there exists a system of local coordinates
$\chi $ near $(t,r,\theta )=(0,0,\theta _{0})$ such that $\Psi
\simeq \chi _{0} \chi _{1}^{2}\chi _{2}\pm \chi _{1}^{k}\chi_{3}$
in a neighborhood of $(\chi _{0},\chi _{1},\chi
_{2},\chi_{3})=(0,0,0,0).$
\end{lemma}
{\it Proof.} Near $\theta _{0}$ we can complete $(Q_{0}(\theta
),R_{k}(\theta ))$ to a system of coordinates on the sphere. Up to
a permutation of the $\theta_i$, we can choose
\begin{equation*}
\begin{array}{l}
(\chi _{0},\chi _{1})(t,r,\theta )=(t,r), \\
\chi _{2}(t,r,\theta )=Q_{t}(\theta )+rh_{1}(t,r,\theta ), \\
\chi _{3}(t,r,\theta )=R_{k}(\theta )+r\tilde{R}(r,\theta ), \\
(\chi _{4},...,\chi _{2n})(t,r,\theta )=(\theta _{3},...,\theta
_{2n-1}).
\end{array}
\end{equation*}
Then the corresponding Jacobian is
\begin{equation*}
|\frac{D\chi }{D(t,r,\theta )}|(0,0,\theta _{0})=\left| \left(
\begin{array}{cc}
\dfrac{\partial \chi _{2}}{\partial \theta _{1}} & \dfrac{\partial \chi _{2}}{%
\partial \theta _{2}} \\
\dfrac{\partial \chi _{3}}{\partial \theta _{1}} & \dfrac{\partial \chi _{3}}{%
\partial \theta _{2}}
\end{array}
\right) \right| (0,0,\theta _{0}) =\left| \left(
\begin{array}{cc}
\dfrac{\partial Q_{0}}{\partial \theta _{1}} & \dfrac{\partial
Q_{0}}{\partial
\theta _{2}} \\
\dfrac{\partial R_{k}}{\partial \theta _{1}} & \dfrac{\partial
R_{k}}{\partial \theta _{2}}
\end{array}
\right) \right|(\theta _{0}) \neq 0.
\end{equation*}
In these new coordinates we have $\Psi (t,z)= (\chi _{0}\chi
_{1}^{2}\chi _{2}+\chi _{1}^{k}\chi _{3})(t,r,\theta ).$
\hfill{$\blacksquare$}
\begin{remark} \label{sur amplitude}
\rm{Coordinates $\chi$ form a system of admissible charts near
$(T,z_{0})$ and these are singular in $z=z_0$ as coordinates on
$T^{\ast}\mathbb{R}^{n}$. In the three systems of coordinates the
measures are $r^{2n-1}|\frac{D\chi }{D(t,r,\theta )}(t,r,\theta
)|dtdrd\theta$, this term $r^{2n-1}$ plays a major role since the
critical sets of our normal forms are $\{r=0\}.$}
\end{remark}
Combining Lemma \ref{Lemme Carte1}, \ref{Lemme Carte2} and
\ref{Lemme Carte3} gives
\begin{theorem}
\label{Formes normales k=entier} If $(H_{3})$ and conditions of
Lemma \ref{Lemme Carte3} are satisfied and if $T$ is a total
period of $d\Phi _{t}(z_{0})$, the phase function $S(t,x,\xi
)-\left\langle x,\xi \right\rangle +tE_{c}$ has one of the
following normal forms on the blow-up of $(T,x_{0},\xi _{0})$ :
\begin{gather}
\text{first normal forms } :(\pm \chi _{0}\chi _{1}^{2})\text{ near }(U\backslash C_{Q_{0}})\text{,} \\
\text{second normal forms } :(\chi _{0}\chi _{1}^{2}\chi _{2}\pm
\chi _{1}^{k})\text{ near }
(C_{Q_{0}}\backslash C_{R_{k}})\text{,} \\
\text{third normal forms } :(\chi _{0}\chi _{1}^{2}\chi _{2}\pm
\chi _{1}^{k}\chi _{3})\text{ near } C_{Q_{0}}\cap
C_{R_{k}}\text{.}
\end{gather}
\end{theorem}
We end this section with two lemmas on asymptotics of oscillatory
integrals.
\begin{lemma} \label{Theo IO 1ere carte}There is a
sequence $(c_{j})_{j}\in \mathcal{D}^{\prime
}(\mathbb{R}^{+}\times \mathbb{R)}$ such that for $k>1$
\begin{equation}
\int\limits_{0}^{\infty }(\int\limits_{\mathbb{R}}e^{i\lambda
tr^{k}}a(t,r)dt)dr \sim \sum\limits_{j=0}^{\infty }\lambda
^{-\frac{j+1}{k}}c_{j}(a),
\end{equation}
where :
\begin{equation}
 c_{l} =\frac{(-1)^{l}}{k}\frac{1}{l!}(\Lambda
_{l}(t)\otimes \delta_{0}^{(l)}(r)),
\end{equation}
with $\Lambda _{l}=\mathcal{F}(x_{-}^{\frac{l+1-k}{k}})$ and
$x_{-}=\max (-x,0)$.
\end{lemma}
{\it Proof.} We define $\hat{g}(\tau
,r)=\mathcal{F}_{t}(a(t,r))(\tau ),$ where $\mathcal{F}_{t}$ is
the partial Fourier transform with respect to $t$. Then we obtain
\begin{equation}
\int\limits_{0}^{\infty }(\int\limits_{\mathbb{R}}e^{i\lambda
tr^{k}}a(t,r)dt)dr=\int\limits_{0}^{\infty }\hat{g}(-\lambda
r^{k},r)dr = \lambda ^{-\frac{1}{k}}\int\limits_{0}^{\infty }\hat{g}%
(-r^{k},\frac{r}{\lambda ^{\frac{1}{k}}})dr.
\end{equation}
Taking the Taylor series in $r$ of $\hat{g}(\tau ,r)$, at the
origin, gives
\begin{equation*}
\hat{g}(-r^{k},\frac{r}{\lambda ^{\frac{1}{k}}})=\sum\limits_{l=0}^{N}\frac{%
\lambda ^{-\frac{l}{k}}}{l!}r^{l}\frac{\partial ^{l}\hat{g}}{\partial r^{l}}%
(-r^{k},0)+\lambda ^{-\frac{N+1}{k}}R(r,\lambda ),
\end{equation*}
Straightforward computations shows that
\begin{equation*}
\lambda ^{-\frac{1}{k}}\int\limits_{0}^{\infty }\hat{g}(-r^{k},\frac{r}{%
\lambda ^{\frac{1}{k}}})dr=\frac{1}{k}\sum\limits_{l=0}^{N}\frac{\lambda ^{-(%
\frac{1+l}{k})}}{l!}\int\limits_{-\infty }^{0}\frac{\partial ^{l}\hat{g}}{%
\partial r^{l}}(r,0)|r|^{\frac{l+1-k}{k}}dr+\mathcal{O}(\lambda ^{-\frac{N+1}{k}}).
\end{equation*}
With $x_{-}=\max (-x,0)$ and $\Lambda _{l}(r)=\mathcal{F}(x_{-}^{%
\frac{l+1-k}{k}})(r)$ the lemma follows. \hfill{$\blacksquare$}
\begin{lemma}
\label{Theo IO 2eme carte}There is a sequence $(c_{j})_{j}\in \mathcal{D}%
^{\prime }(\mathbb{R}^{+}\times \mathbb{R)}$ such that for $k>1$
\begin{equation}
\int\limits_{0}^{\infty }e^{i\lambda r^{k}}a(r)dr \sim
\sum\limits_{j=0}^{\infty }\lambda ^{-\frac{j+1}{k}}c_{j}(a),
\end{equation}
with :
\begin{equation}
c_{j}(a) =\frac{(-1)^{j}}{k}\Gamma (\frac{j}{k})\exp (i\pi \frac{j}{2k})%
\frac{\partial ^{j}a}{\partial r^{j}}(0).
\end{equation}
\end{lemma}
{\it Proof.} We use the Berstein-Sato polynomial, see e.g.
\cite{WON}. We write
\begin{equation*}
\int\limits_{0}^{\infty }e^{i\lambda r^{k}}a(r)dr=\frac{1}{2i\pi }%
\int\limits_{\gamma }e^{i\pi \frac{z}{2}}\Gamma (z)\lambda
^{-z}(\int\limits_{0}^{\infty }a(r)r^{-kz}dr)dz,
\end{equation*}
with $\gamma =]c-i\infty ,c+i\infty \lbrack ,$
$\text{\textrm{Re}}(c)<k^{-1}.$ Since for all positive $r$ we have
\begin{equation*}
\frac{\partial ^{k}}{\partial r^{k}}(r^{k})^{1-z}=r^{-kz}\prod
\limits_{j=1}^{k}(j-kz),
\end{equation*}
we can compute the asymptotic by the residue method. All poles are
simple and, by pushing of the complex path of integration to the
right, we obtain
\begin{equation*}
\lim_{z\rightarrow \frac{l}{k}}\frac{(z-\frac{l}{k})}{\prod
\limits_{j=1}^{k}(j-kz)}e^{i\pi \frac{z}{2}}\Gamma (z)\lambda
^{-z}(-1)^{k}\int\limits_{0}^{\infty }\frac{\partial
^{k}a}{\partial r^{k}}(r)r^{k-kz}dr = \mu _{l}\frac{\partial
^{l-1}a}{\partial r^{l-1}}(0).
\end{equation*}
Straightforward computations then show that $\mu
_{l}=(-1)^{k}\Gamma(\frac{l}{k})\exp (i\pi \frac{l}{2k})$.
\hfill{$\blacksquare$}
\section{Proofs of the main theorems}
We start with the simpler case of $T$ being a total period of
$d\Phi _{t}(z_{0})$. Afterwards we study the contributions of
non-total periods $T$, distinguishing out particular case of a
function $p_{2}$ whose restriction to the linear subspace
$\mathbb{\frak{F}}_{T}$ has constant sign. In the following we
suppose, without loss of generality, that the support of the
amplitude contains only one non-zero period of the linearized
flow.
\subsection{Blow-up and partition of the sphere}
We will apply the results of sections 4 and 5, and we recall that
\begin{equation*}
R(x,\xi ) =S(T,x,\xi )-\left\langle x,\xi \right\rangle
,\end{equation*}
\begin{equation*}
(t-T)g(t,x,\xi )+tE_{c} =(t-T)(Q_{t}(x,\xi )+h(t,x,\xi )).
\end{equation*}
By a time translation and with the polar coordinates $(x,\xi
)=r\theta$, $\theta \in \mathbb{S}^{2n-1}(\mathbb{R)}$, we obtain
for the top order part of $(2\pi h)^{n}\gamma_{2}(E_{c},h)$ :
\begin{equation*}
I(T,h)=\int\limits_{\mathbb{R\times \lbrack }0,\infty \lbrack }\int\limits_{%
\mathbb{S}^{2n-1}}a(t,r\theta )e^{\frac{i}{h}\Psi (t,r,\theta
)}r^{2n-1}dtdrd\theta ,
\end{equation*}
where $d\theta $ is the standard surface measure on the
sphere.\medskip \newline%
\textbf{Partition of unity on the sphere.}\medskip \newline%
Since $C_{Q_{T}}\cap\mathbb{S}^{2n-1}(\mathbb{R)}$ is compact, we
can introduce a finite partition of unity
\begin{equation*}
\sum_{i\in I}\Omega _{i}^{1}(\theta )+\sum_{j\in J}\Omega
_{j}^{2}(\theta )+\sum_{l\in L}\Omega _{l}^{3}(\theta )=0\text{ on
}\mathbb{S}^{2n-1}(\mathbb{R)}\text{,}
\end{equation*}
with the property that $Q_{T}(\theta ) \neq 0 \text{ on}\bigcup
\rm{supp} (\Omega _{i}^{1})$, $R_{k}(\theta ) \neq 0\text{ on
}\bigcup \rm{supp} (\Omega _{j}^{2})$ and $C_{Q_{T}}\cap C_{R_{k}}
\subset \bigcup \rm{supp} (\Omega _{l}^{3})$. We split up the
integral $I(T,h)$ according to this partition of unity and use the
normal forms of Theorem \ref{Formes normales k=entier}. On any
chart let be $J\chi $ the Jacobian of the relevant diffeomorphism
of blow-up from Theorem \ref{Formes normales k=entier}. For $j$ in
$I$, $J$ and $L$ respectively, we define
\begin{equation}
a_{i,j}(\chi _{0},\chi _{1},...,\chi _{2n})=(\chi ^{-1})^{\ast
}(\Omega _{j}^{i}(\theta )a(t,r\theta )r^{2n-1}|J\chi
|^{-1}),\text{ } i\in\{1,2,3\},
\end{equation}
Then, from Theorem \ref{Formes normales k=entier}, we obtain the
local contributions
\begin{gather*}
\int\limits_{\mathbb{R\times R}^{+}\times \mathbb{S}^{2n-1}}\Omega
_{i}^{1}(\theta )a(t,r\theta )e^{\frac{i}{h}\Psi (t,r\theta
)}r^{2n-1}dtdrd\theta \\
=\int\limits_{\mathbb{R\times \lbrack
}0,\infty \lbrack }e^{\pm\frac{i}{h}\chi _{0}\chi
_{1}^{2}}A_{1,i}(\chi _{0},\chi _{1})d\chi _{0}d\chi _{1},
\end{gather*}
for the first normal forms. Also
\begin{equation*}
 \int\limits_{\mathbb{R\times
R}^{+}\times \mathbb{S}^{2n-1}}\Omega _{j}^{2}(\theta )a(t,r\theta
)e^{\frac{i}{h}\Psi (t,r\theta )}r^{2n-1}dtdrd\theta
\end{equation*}
\begin{equation*}
=\int\limits_{\mathbb{R\times R}^{+}\times
\mathbb{R}}A_{2,j}(\chi_{0},\chi_{1},\chi_{2})e^{\frac{i}{h}(\chi_{0}\chi
_{1}^{2}\chi _{2}\pm \chi _{1}^{k})}d\chi _{0}d\chi _{1}d\chi
_{2},
\end{equation*}
for the second normal forms and for the third normal forms
\begin{equation*}
\int\limits_{\mathbb{R\times R}^{+}\times \mathbb{S}^{2n-1}}\Omega
_{l}^{3}(\theta )a(t,r\theta )e^{\frac{i}{h}\Psi (t,r\theta
)}r^{2n-1}dtdrd\theta
\end{equation*}
\begin{equation*}
=\int\limits_{\mathbb{R\times \mathbb{R}^{+}}\times
\mathbb{R}^{2}}A_{3,l}(\chi _{0},..,\chi _{3})e^{\frac{i}{h}(\chi
_{0}\chi _{1}^{2}\chi _{2}\pm \chi _{1}^{k}\chi _{3})}d\chi
_{0}..d\chi _{3},
\end{equation*}
where the amplitudes are respectively given by
\begin{equation}
A_{1,i}(\chi _{0},\chi _{1})=\int a_{2,i}(\chi _{0},\chi
_{1},...,\chi _{2n})d\chi _{2}...d\chi _{2n},  \label{Amplitude1}
\end{equation}
\begin{equation}
A_{2,j}(\chi _{0},\chi _{1},\chi _{2})=\int a_{2,j}(\chi _{0},\chi
_{1},...,\chi _{2n})d\chi _{3}...d\chi _{2n}, \label{Amplitude2}
\end{equation}
\begin{equation}
A_{3,l}(\chi _{0},\chi _{1},\chi _{2},\chi _{3})=\int a_{3,l}(\chi
_{0},\chi _{1},...,\chi _{2n})d\chi _{4}...d\chi _{2n}.
\label{Amplitude3}
\end{equation}
It is convenient, for the calculations below, to introduce
$A_{j,i}=\chi _{1}^{2n-1}\tilde{A}_{j,i}$, for $j\in \{1,2,3\}$,
c.f. Remark [\ref{sur amplitude}]. By construction the functions
$A_{j,i}$ are of compact support in their system of coordinates.
\begin{remark}
\rm{We obtain for each new phases the following critical sets}
\begin{equation*}
\left\{
\begin{array}{l}
\frak{C}(\chi _{0}\chi _{1}^{2})=\{\chi _{1}=0\}, \\
\frak{C}(\chi _{0}\chi _{1}^{2}\chi _{2}+\chi _{1}^{k})=\{\chi _{1}=0\}, \\
\frak{C}(\chi _{0}\chi _{1}^{2}\chi _{2}+\chi _{1}^{k}\chi
_{3})=\{\chi _{1}=0\}.
\end{array}
\right.
\end{equation*}
\rm{Where $\frak{C}(f)$ denotes the critical set of a function
$f$.}
\end{remark}
\subsection{Analysis in the case of a total period}
\noindent \textbf{First normal forms.}\newline We note
$\mathcal{F}$ and $\mathcal{F}_{t}$ the total and partial Fourier
transform with respect to $t$. For an amplitude of the form
$a(r,t)=\mathcal{O}(r^{2n-1})$ and $k=2$, where $r=\chi_1$, Lemma
\ref{Theo IO 1ere carte} shows that the first non-zero coefficient
is obtained for $l_{0}=2n-1$. With $\lambda =h^{-1}$, this gives
the contribution
\begin{equation*}
\frac{\lambda ^{-n}}{(2n-1)!}c_{2n-1}(A_{1,i})=\lambda
^{-n}\left\langle \Lambda
_{2n-1}(t),\tilde{A}_{1,i}(t,0)\right\rangle ,
\end{equation*}
where $\Lambda _{2n-1}(t)=\mathcal{F}(x_{-}^{n-1})(t).$ If we
define the distribution
\begin{equation*}
\Lambda_{1,i} (a)=\langle \Lambda _{2n-1}(\chi _{0}),\int(\chi
^{-1})^{\ast }(\Omega _{i}^{1}(\theta )a(t,r\theta )|J\chi
|^{-1})(\chi _{0},0,\chi _{2},...,\chi _{2n})d\chi _{2}...d\chi
_{2n}\rangle ,
\end{equation*}
we obtain, as $\lambda \rightarrow +\infty$, the asymptotic
equivalent
\begin{equation}
\int\limits_{\mathbb{R}}\int\limits_{0}^{\infty }e^{i\lambda
\chi_{0}\chi_{1}^{2}}A_{1,i}(\chi_0,\chi_1)d\chi_{0}d\chi_{1}=
-\frac{1}{2}\Lambda_{1,i} (a)\lambda ^{-n}+\mathcal{O}(\lambda ^{-n-\frac{1%
}{2}}).  \label{DA sur carte 1}
\end{equation}
\textbf{Second normal forms.}\newline If we write $a$ for
$\tilde{A}_{2,j}$, $\lambda=h^{-1}$ and $(t,r,v)$ for
$(\chi_0,\chi_1,\chi_2)$, then we have to analyze the asymptotic
behavior of the oscillatory integral
\begin{equation}
I_{2}(\lambda)=\int\limits_{0}^{\infty
}\int\limits_{\mathbb{R}^{2}}e^{i\lambda
(tr^{2}v+r^{k})}a(t,r,v)dtdvr^{2n-1}dr
\end{equation}
If we let $\hat{a}(\tau ,r,v)=\mathcal{F}_{t}(a(t,r,v))(\tau )$,
the Fourier transform with respect to $t$, then we find by easy
manipulations, that
\begin{equation}\label{Equationrepport}
I_{2}(\lambda)=\lambda ^{-1}\int\limits_{0}^{\infty} r^{2n-3}
e^{i\lambda
r^{k}}\int\limits_{\mathbb{R}}\hat{a}(-w,r,\frac{w}{\lambda
r^{2}})dwdr,
\end{equation}
where we made the change of variables $(v,r)\rightarrow (w,r)$,
$w=\lambda r^2 v$.\\
We now Taylor expand up till order $n-3$ :
\begin{equation*}
\hat{a}(-w,r,\dfrac{w}{\lambda r^{2}})=\sum\limits_{j=0}^{n-3}c_{j}(-w,r)(%
\dfrac{w}{\lambda r^{2}})^{j}+R_{n-2}(r,w,\lambda ),
\end{equation*}
where
\begin{equation*}
c_{j}(w,r)=\dfrac{1}{j!}\dfrac{\partial ^{j}\hat{a}}{\partial
v^{j}}(w,r,0),
\end{equation*}
and
\begin{equation*}
R_{n-2}(r,w,\lambda )=(\frac{w}{\lambda
r^2})^{n-2}\int\limits_{0}^{1}
\partial_{r}^{n-2}\hat{a}(-w,r,\frac{tw}{\lambda
r^2})\frac{(1-t)^{n-3}}{(n-3)!}dt.
\end{equation*}
Substitution in Eq.(\ref{Equationrepport}) leads to
\begin{equation*}
I_2(\lambda)=\sum\limits_{j=0}^{n-3} J_{j}(\lambda
)+r_{n-2}(\lambda),
\end{equation*}
where
\begin{gather*}
J_{j}(\lambda)=\lambda^{-(1+j)} \int\limits_{0}^{\infty}
e^{i\lambda r^{k}}a_{j}(r)r^{2n-3-2j}dr, \\
a_{j}(r)=\int\limits_{\mathbb{R}}c_{j}(-w,r)w^{j}dw=\frac{(-i)^{j}}{j!}
((\dfrac{\partial ^{2}}{\partial v\partial t})^{j}a)(0,r,0),
\end{gather*}
and
\begin{equation*}
r_{n-2}(\lambda)=\lambda^{-1}\int\limits_{0}^{\infty}\int\limits_{\mathbb{R}}
R_{n-2}(r,w,\lambda)e^{i\lambda r^{k}}r^{2n-3}drdw.
\end{equation*}
Observe that the oscillatory integrals $J_j(\lambda)$ can be
treated individually using Lemma  \ref{Theo IO 2eme carte}, but we
first analyze the remainder term, $r_{n-2}(\lambda)$, which equals
\begin{equation}\label{Equationrepport1}
\frac{1}{\lambda^{n-1}}\int\limits_{0}^{\infty}re^{i\lambda r^{k}}
\int\limits_{\mathbb{R}}\int\limits_{0}^{1}w^{n-2}\partial_{v}^{n-2}
\hat{a}(-w,r,\frac{tw}{\lambda
r^2})\frac{(1-t)^{n-3}}{(n-3)!}dtdwdr.
\end{equation}
We split the integral with respect to $dr$ as :
\begin{equation*}
\int\limits_{0}^{\infty} dr=\int\limits_{0}^{A}
dr+\int\limits_{A}^{\infty} dr,
\end{equation*}
where $A=A(\lambda)$ will be chosen below. We accordingly split
$r_{n-2}(\lambda)$ as
$r_{n-2}^{1,A}(\lambda)+r_{n-2}^{2,A}(\lambda)$ and
$r_{n-2}^{1,A}(\lambda)$ is given by Eq.(\ref{Equationrepport1}),
with the integral restrained to $[0,A]$. Easy estimates then show
that $r_{n-2}^{1,A}(\lambda) \leq C A^2 \lambda^{-(n-1)}$ with $C$
independent of $\lambda$.

Next, for $r_{n-2}^{2,A}(\lambda)$ we do an integration by part
with respect to $t$, this leads to
\begin{gather}
r_{n-2}^{2,A}(\lambda)=\lambda^{-(n-1)}\int\limits_{A}^{\infty} r
e^{i\lambda r^{k}}\int\limits_{\mathbb{R}}
\frac{w^{n-2}}{(n-2)!}\partial_{v}^{n-2}\hat{a}(-w,r,0)dwdr \notag \\
+\lambda^{-n}\int\limits_{A}^{\infty}\frac{e^{i\lambda r^{k}}}{r}
\int\limits_{\mathbb{R}}\int\limits_{0}^{1}w^{n-1}\partial_{v}^{n-1}
\hat{a}(-w,r,\frac{tw}{\lambda r^2})
\frac{(1-t)^{n-2}}{(n-2)!}dtdwdr. \label{Equationrepport2}
\end{gather}
We then observe that the first integral in
Eq.(\ref{Equationrepport2}) is equal to
\begin{gather*}
\lambda^{-(n-1)}\int\limits_{0}^{\infty} r e^{i\lambda
r^{k}}\int\limits_{\mathbb{R}}\frac{w^{n-2}}{(n-2)!}\partial_{v}^{n-2}\hat{a}(-w,r,0)dwdr
+\mathcal{O}(A^2 \lambda^{-(n-1)})\\
=J_{n-2}(\lambda)+\mathcal{O}(A^2 \lambda^{-(n-1)}),
\end{gather*}
by similar estimates as for $r_{n-2}^{1,A}(\lambda)$. Finally, the
last integral in Eq.(\ref{Equationrepport2}) can be estimated by
\begin{equation*}
\frac{\lambda^{-n}}{(n-2)!}\int\limits_{A}^{\infty}\int\limits_{\mathbb{R}}
\frac{|w|^{n-1}}{r}||\partial_{v}^{n-1}\hat{a}(-w,r,\bullet)||_{\infty}dwdr\leq
C \lambda^{-n} |\mathrm{log}(A)|,
\end{equation*}
remembering that $\hat{a}$ has a compact support in $v$. In
conclusion, we find that
\begin{equation*}
I_2(\lambda)=\sum\limits_{j=0}^{n-2} J_j(\lambda)
+\mathcal{O}(\lambda^{-(n-1)}A^2)+\mathcal{O}(\lambda^{-n}
\mathrm{log} (A))=\sum\limits_{j=0}^{n-2}
J_j(\lambda)+\mathcal{O}(\lambda^{-n} \mathrm{log} (\lambda)),
\end{equation*}
where we have chosen $A=\lambda^{-\frac{1}{2}}$.

Now by Lemma \ref{Theo IO 1ere carte}, and since
$a_0(r)=\tilde{A}_{2,j}(0,r,0)$, we obtain
\begin{equation}
J_0(\lambda)=\mu _{0}\lambda ^{\frac{2-k-2n}{k}}\frac{\partial
^{2n-3}a_{0}}{\partial r^{2n-3}}(0)
+\mathcal{O}(\lambda^{\frac{1-2n-k}{k}}),
\end{equation}
with $\mu _{0}=-\frac{1}{k} \Gamma
(\frac{n-2}{k})\mathrm{exp}(i\pi\frac{n-1}{k})$. While, if $j>0$,
the same lemma shows that
\begin{equation*}
J_j(\lambda)=\mathcal{O}(\lambda
^{-\frac{2n+(j+1)(k-2)}{k}}),\text{ } 0<j\leq n-2.
\end{equation*}
All of the latter are dominated by the top order term of
$J_0(\lambda)$, and dominate the remainder, which is
$\mathcal{O}(\lambda^{-n}\mathrm{log}(\lambda))$, since for
$j=n-2$ and for all $k\geq 3$ :
\begin{equation*}
\frac{2n+(k-2)(j+1)}{k}=n-1+\frac{2}{k}<n.
\end{equation*}
As concerns the error term in (\ref{Equationrepport2}), this also
dominates the remainder term, since $n>\frac{k-2}{k-1}$ for all
$n\geq 2$. Since $h=\lambda^{-1}$, we have shown that the
contribution of the second normal form to the asymptotics of
$I(T,h)$ is :
\begin{equation}
I_2(h)=\mu_0 \tilde{A}_{2,j}(0,0,0)
h^{\frac{2n+k-2}{k}}+\mathcal{O}(h^{\frac{2n+k-1}{k}}).
\end{equation}
\textbf{Third normal forms.}\newline%
We use the same strategy as for second normal forms. If
$\lambda=h^{-1}$, $(t,r,v,s)=(\chi_0,\chi_1,\chi_2,\chi_3)$ and
$a=\tilde{A}_{3,j}$, we write
\begin{equation*}
I_3(\lambda)=\int\limits_{0}^{\infty
}(\int\limits_{\mathbb{R}^{3}}e^{i\lambda
(tr^{2}v+r^{k}s)}a(t,r,v,s)dtdvds)r^{2n-1}dr,
\end{equation*}
with $\hat{a}(\tau ,r,v,s)=\mathcal{F}_{t}(a(t,r,v,s))(\tau )$, a
Taylor expansion, up till order $n-3$, gives again
\begin{equation*}
\hat{a}(-w,r,\dfrac{w}{\lambda
r^{2}},s)=\sum\limits_{j=0}^{n-3}c_{j}(-w,r,s) (\dfrac{w}{\lambda
r^{2}})^{j}+R_{n-2}(r,w,\lambda ,s),
\end{equation*}
with
\begin{gather*}
R_{n-2}(r,w,\lambda ,s)=(\frac{w}{\lambda
r^2})^{n-2}\int\limits_{0}^{1}
\partial_{r}^{n-2}\hat{a}(-w,r,\frac{tw}{\lambda r^2},s)\frac{(1-t)^{n-3}}{(n-3)!}dt,\\
c_{j}(w,r,s)=\dfrac{1}{j!}\dfrac{\partial ^{j}\hat{a}}{\partial
v^{j}}(w,r,0,s),\\
a_{j}(r,s)=\int\limits_{\mathbb{R}}c_{j}(-w,r,s)w^{j}dw=(-i)^{j}
((\dfrac{\partial ^{2}}{\partial v\partial t})^{j}a)(0,r,0,s).
\end{gather*}
This leads to
\begin{gather*}
I_3(\lambda)= \sum\limits_{j=0}^{n-3}
K_j(\lambda)+r_{n-2}(\lambda )\\
=\sum\limits_{j=0}^{n-3}\lambda ^{-(1+j)}\int e^{i\lambda
r^{k}s}a_{j}(r,s)r^{2n-3-2j}drds+r_{n-2}(\lambda ).
\end{gather*}
The $K_j(\lambda)$ can be treated individually via Lemma \ref{Theo
IO 1ere carte}. For $j=0$, we obtain
\begin{gather*}
K_0(\lambda)=\lambda ^{-1}\int\limits_{0}^{\infty
}\int\limits_{\mathbb{R}}e^{i\lambda
r^{k}s}a_{0}(r,s)r^{2n-3}drds%
\\ \sim \lambda ^{-1}(\sum_{l=0}^{N}\lambda
^{-\frac{l+1}{k}}c_{l}(a_{0}(r,s)r^{2n-3})+\mathcal{O}(\lambda ^{-\frac{N+2}{k}%
})) .
\end{gather*}
The leading term, obtained for $l=2n-3$, is
\begin{equation*}
\lambda ^{-(\frac{2n-2+k}{k})}c_{2n-3}(a_{0}(r,s)r^{2n-3})=-\frac{1}{k}%
\lambda ^{-(\frac{2n-2+k}{k})}\left\langle \mathcal{F}(x_{-}^{\frac{2n-2-k}{k%
}})(s),a_{0}(0,s)\right\rangle .
\end{equation*}
Hence for our amplitude we have the main contribution
\begin{equation}
-\frac{1}{k}\lambda ^{-(\frac{2n-2+k}{k})}\left\langle \mathcal{F}(x_{-}^{%
\frac{2n-2-k}{k}})(\chi _{3}),\tilde{A}_{3,l}(0,0,0,\chi
_{3})\right\rangle +\mathcal{O}(\lambda ^{-(\frac{2n-1+k}{k})}).
\label{contribution carte 3}
\end{equation}
Like for the second normal forms, the other terms $K_j(\lambda)$,
with $j>0$, and the remainder $r_{n-2}(\lambda)$ give
contributions of strictly lower orders.

Finally, on each local chart the main contributions are
\begin{equation*}
\left\{
\begin{array}{l}
\text{first normal forms }:c_{0,1}(a)\lambda ^{-n}, \\
\text{second normal forms }:c_{0,2}(a)\lambda ^{\frac{2-2n-k}{k}}, \\
\text{third normal forms }:c_{0,3}(a)\lambda ^{\frac{2-2n-k}{k}}.
\end{array}
\right.
\end{equation*}
The contributions of charts 2 and 3 are dominant, since
$\frac{2-2n-k}{k}>-n$, $\forall k>2$.
\begin{remark}
\rm{The proofs above show in fact much more than just an
asymptotic equivalent for $I(T,h)$, and therefore for $\gamma_2
(E_c ,h)$. They show the existence of a limited asymptotic
expansion in the case of indefinite $Q_T$ :
\begin{equation}
I(T,h)=\sum c_{\nu} h^{\frac{2n+k-2+\nu}{k}}+\mathcal{O}(h^n
|\mathrm{log}(h)|),
\end{equation}
where the sum is over all $\nu$ such that
$\frac{2n+k-2+\nu}{k}<n$, or $\nu<(k-2)(n-1)$, and of a complete
asymptotic expansion if $Q_T$ is definite. A similar remark
applies for the case of a non-total period, which we examine in
the next section.}
\end{remark}

We now compute the leading term of the expansion in case of a
non-definite $Q_T$ and for $T$ a total period of the linearized
flow.\medskip\newline %
\textbf{Case of an empty intersection of cones.}\\
With $R_{k}\neq 0$ on $C_{p_{2}}\backslash\{0\}$, we can assume
that $R_{k}$ is positive on $C_{p_{2}}\backslash \{0\}$. The main
contribution is here given by the second normal form and is
\begin{gather*}
\int\limits_{\mathbb{R\times \lbrack }0,\infty \lbrack \times \mathbb{R}%
}A_{2,j}(\chi _{0},\chi _{1},\chi _{2})e^{i\lambda (\chi _{0}\chi
_{1}^{2}\chi _{2}+\chi _{1}^{k})}d\chi _{0}d\chi _{1}d\chi
_{2}\\
=\mu_{k}\lambda ^{\frac{2-k-2n}{k}}\tilde{A}_{2,j}(0,0,0)
+\mathcal{O}(\lambda ^{\frac{2-k-2n}{k}-\frac{1}{k}}),
\end{gather*}
with $\mu_{k}$ given by Lemma \ref{Theo IO 2eme carte}. By
definition of $\chi$, the amplitude $\tilde{A}_{2,j}(0,0,0)$ is
\begin{equation*}
\int(\chi ^{-1})^{\ast }(\Omega _
{j}^{2}(\theta )|R_{k}(\theta )+r\tilde{R}(r,\theta )|^{-\frac{(2n-1)}{k}%
}a(t+T,r\theta )|J\chi |^{-1})(0,\chi _{3},..,\chi _{2n})d\chi
_{3}..d\chi _{2n}.
\end{equation*}
For $z\in\mathbb{R}^3$, we write this delta-Dirac distribution as
an oscillatory integral
\begin{equation*}
\frac{1}{(2\pi )^{3}}\int e^{-i\left\langle z,(\chi _{0},\chi
_{1},\chi _{2})(t,r,\theta )\right\rangle }\Omega _
{j}^{2}(\theta )|R_{k}(\theta )+r\tilde{R}(r,\theta )|^{-\frac{2n-1}{k}}%
a(t+T,r\theta )dzdtdrd\theta .
\end{equation*}
If we use $y_{2}=z_{2}|R_{k}(\theta )+r\tilde{R}(r,\theta )|^{%
\frac{1}{k}}$, integration in $(y_{2},r)$ gives
\begin{equation*}
(2\pi )^{2}\tilde{A}_{2,j}(0,0,0)=\int e^{-i\left\langle
(z_{1},z_{3}),(\chi _{0},\chi _{2})(t,0,\theta )\right\rangle
}\Omega _{j}^{2}(\theta )|R_{k}(\theta
)|^{-\frac{2n}{k}}a(t+T,0)dtd\theta dz_{1}dz_{3}.
\end{equation*}
Since $\chi _{0}(t,0,\theta )=t|R_{k}(\theta )|^{-\frac{2}{k}}$,
with $y_{1}=z_{1}|R_{k}(\theta )|^{-\frac{2}{k}}$, by integration
in $(t,y_{1})$
\begin{equation*}
(2\pi )\tilde{A}_{2,j}(0,0,0)=a(T,0)\int e^{-iz_{3}p_{2}(\theta
)}\Omega _{j}^{2}(\theta )|R_{k}(\theta
)|^{-\frac{2n-2}{k}}d\theta dz_{3}.
\end{equation*}
We define the Liouville measure $dL_{p_{2}}$ on $C_{p_{2}}\cap
\mathbb{S}^{2n-1}$ via $dL_{p_{2}}(\theta )\wedge dp_{2}(\theta
)=d\theta$. Since charts associated to second normal forms cover
the trace of the cone, by summation over the partition of unity we
obtain
\begin{equation*}
I(T,h)=h^{\frac{2n+k-2}{k}}(\mu_{k}
a(T,0)\int\limits_{C_{p_{2}}\cap \mathbb{S}^{2n-1}}|R_{k}(\theta
)|^{-\frac{2n-2}{k}}dL_{p_{2}}(\theta )
+\mathcal{O}(h^{\frac{1}{k}})),
\end{equation*}
where $\theta $ are now local coordinates on the surface $%
C_{p_{2}}\cap \mathbb{S}^{2n-1}$. The top-order contribution to
the trace formula follows from $\gamma_{2}(E_{c},T,h)=(2\pi
)^{-n-1}h^{-n}I(T,h)$. \medskip\newline%
\textbf{Case of a non-empty intersection of cones.}\newline%
Here the main contribution is given by the normal forms 2 and 3.
The local contribution of any chart associated to second normal
forms can be computed like in the previous section. The
contribution of the third normal forms is given by equation
(\ref{contribution carte 3}) with the amplitude
$\tilde{A}_{3,l}(0,0,0,\chi _{3})$. Let be
$z=(z_{1},z_{2},z_{3})$, we use again an oscillatory
representation of the delta-Dirac distribution via
\begin{equation*}
\frac{1}{(2\pi )^{3}}\int e^{-i\left\langle z,(\chi _{0},\chi
_{1},\chi _{2})\right\rangle }\left\langle
\mathcal{F}(x_{-}^{\frac{2n-2-k}{k}})(\chi
_{3}),\tilde{A}_{3,l}(\chi _{0},..,\chi _{3})\right\rangle d\chi
_{0}d\chi _{1}d\chi _{2}dz.
\end{equation*}
Since $(\chi _{0},\chi _{1})=(t,r)$, integration w.r.t.
$(z_{1},z_{3},\chi _{0},\chi _{1})$ gives
\begin{equation*}
2\pi \left\langle \mathcal{F}(x_{-}^{\frac{2n-2-k}{k}})(\chi _{3}),\tilde{A}%
_{3,l}(0,0,0,\chi _{3})\right\rangle
\end{equation*}
\begin{equation*}
=\int e^{-iz_{2}\chi _{2}(0,0,\theta )}\left\langle
\mathcal{F}(x_{-}^{\frac{2n-2-k}{k}})(\chi _{3}(0,0,\theta
)),\Omega _{l}^{3}(\theta )a(T,0)\right\rangle d\theta dz_{2}.
\end{equation*}
A classical result, see \cite{HOR1} volume 1 page 167, is
\begin{equation}
\mathcal{F}(x_{-}^{\frac{2n-2-k}{k}})(\chi _{3})=\Gamma
(\frac{2n-2}{k})\exp (i\pi \frac{n-1}{k})(\chi
_{3}+i0)^{-\frac{2n-2}{k}}.
\end{equation}
By construction $(\chi _{2},\chi _{3})(0,0,\theta )=(p_{2}(\theta
),R_{k}(\theta ))$ and with the notations of Theorem \ref{Main
result Q non positive} we obtain, using again $p_{2}$ as a local
coordinate
\begin{equation*}
\mu _{k}a(T,0)\int\limits_{C_{p_{2}}\mathbb{\cap
S}^{2n-1}}(R_{k}(\theta)+i0)^{-\frac{2n-2}{k}}\Omega_{l}^{3}(\theta
)dL_{p_{2}}(\theta ).
\end{equation*}
But $(R_{k}(\theta )+i0)^{-\frac{2n-2}{k}}=R_{k}(\theta )^{-\frac{2n-2}{k}%
} $ on a chart associated to the second normal form and the total
contribution arising from normal forms 2 and 3 is
\begin{equation}
I(T,h)=h^{\frac{2n-2+k}{k}}(\mu
_{k}a(T,0)\int\limits_{C_{p_{2}}\mathbb{\cap
S}^{2n-1}}(R_{k}(\theta )+i0)^{-\frac{2n-2}{k}}dL_{p_{2}}(\theta )
+\mathcal{O}(h^{\frac{1}{k}})).
\end{equation}
Since $a(T,0)=\hat{\varphi}(T)\exp (iTp^{1}(z_{0}))$, see formula
(\ref{density}) below, this proves Theorem \ref{Main result Q non
positive} for a total period. \hfill{$\blacksquare$}
\subsection{Case of a non-total period}
Let be $T$ a non total period of the linearized flow. We can
assume, up to a permutation of coordinates, that
$\mathbb{\frak{F}}_{T}=\{z=(z_{1},z_{2})\in
\mathbb{R}^{l_{T}}\times \mathbb{R}^{2n-l_{T}}\text{ }/\text{
}z_{2}=0\}$. We can apply the Morse lemma with parameter since the
phase function at time $T$ is only degenerate in $z$ along
$\mathbb{\frak{F}}_{T}$ (cf. Corollary \ref{theo periode flot
linéarise}). The quadratic part $S_{2}(t,x,\xi )$ of the function
$S$ is given by
\begin{equation*}
\Phi _{t}(\partial _{\xi }S,\xi )=d\Phi _{t}(0)(\partial _{\xi
}S_{2},\xi )+\mathcal{O}(||(x,\xi )||^{2})=(x,\partial
_{x}S_{2})+\mathcal{O}(||(x,\xi )||^{2}).
\end{equation*}
Using Theorem 5.3 of \cite{HOR2} we can assume that $\det [\frac{%
\partial ^{2}S}{\partial x\partial \eta }]\neq 0$ and the
function $S_{2}$ is well determined locally. Then we have the
following facts
\begin{equation*}
\left\{
\begin{array}{l}
d_{z}\Psi (t,z)=0\Leftrightarrow z=z_{0},\text{ }\forall t, \\
\mathrm{Hess}_{z_{2}}(\Psi )(T,z_{0})= \partial_{z_{2}}
^{2}\Psi(T,z_{0}) \text{, is invertible.}
\end{array}
\right.
\end{equation*}
By the Morse lemma, after a change of variable $z\rightarrow
\tilde{z}$ and calling $\tilde{z}$ again $z$, we have
\begin{equation*}
\Psi (t,z)=q(z_{2})+\Psi (t,z_{1},z_{2}(t,z_{1}))=q(z_{2})+\tilde{\Psi}%
(t,z_{1}),
\end{equation*}
and by Corollary \ref{theo periode flot linéarise} again, $q=\frac{1}{2}q_{T}=p_{2}|%
\mathbb{\frak{F}}_{T}^{\perp }.$ In the following, we note
\begin{gather*}
R(t,z_{1}) \equiv R(z_{1},z_{2}(t,z_{1})),\\
g(t,z_{1}) \equiv g(t,z_{1},z_{2}(t,z_{1})), \\
\Psi (t,z_{1}) \equiv \Psi
(t,z_{1},z_{2}(t,z_{1}))=R(z_{1},z_{2}(t,z_{1}))+(t-T)g(t,z_{1},z_{2}(t,z_{1})).
\end{gather*}
With these conventions we can write
\begin{equation*}
I(T,h)=\int\limits_{z_{1}\in \mathbb{R}^{l_{T}}}\int\limits_{z_{2}\in \mathbb{R}%
^{2n-l_{T}}}e^{\frac{i}{2h}q_{T}(z_{2})}\tilde{a}(t,z_{1},z_{2})e^{\frac{i}{h}%
(R(t,z_{1})+(t-T)g(t,z_{1}))}dtdz_{1}dz_{2},
\end{equation*}
where $\tilde{a}$ is the new amplitude after change of variables
due to the Morse lemma. The stationary phase method applied to the
$z_{2}$-integral gives
\begin{equation}
\int\limits_{z_{2}\in
\mathbb{R}^{2n-l_{T}}}e^{\frac{i}{2h}q_{T}(z_{2})}\tilde{a}
(t,z_{1},z_{2})dz_{2}=\sum\limits_{\nu=0}^{N}c_{\nu}h^{n+\nu-d_{T}}A_{\nu}(t,z_{1})
+\mathcal{O}(h^{N+\nu -d_{T}+1}), \label{Report8}
\end{equation}
and in particular for the leading term we have
\begin{gather*}
A_{0}(t,z_{1})
=a(t,z_{1},0)|\frac{\partial{z_{2}}(t,z_{1})}{\partial
(z_{1},t)}|^{-1},\\
c_{0} =\frac{(2\pi )^{n-d_{T}}\exp
(i\frac{\pi}{4}\mathrm{sgn}(q_{T}))}{|\det(q_{T})|^{\frac{1}{2}}}.
\end{gather*}
We now distinguish two different cases : the quadratic form
$Q_{T}$, that is $p_{2}$ restricted to $\mathbb{\frak{F}}_{T},$ is
definite or non-definite. In the first case only the normal forms
$\pm \chi _{0}\chi _{1}^{2}$ will occur, i.e. the knowledge of
$d\Phi_{t}(z_0)$ is sufficient. In the second case the main
contribution involves $R_k$, and hence the operator
$d^{k-1}\Phi_{t}(z_{0})$. \medskip \newline%
\textbf{The case $Q_{T}$ definite.}\newline Up to a complex
conjugation we can assume that $Q_{T}$ is positive. With polar
coordinates $z_1=(r\theta)$ and $a_{|z_{2}=0}=a(t,r\theta,0)$, the
new amplitude is
\begin{equation}
A_{0}(\chi _{0},\chi _{1})=\int(\chi ^{-1})^{\ast }(a(t,r\theta
,0)r^{l_{T}-1}|J\chi (t,r,\theta )|)d\chi _{2}...d\chi _{2d_{T}},
\end{equation}
with $A_{0}(\chi _{0},\chi _{1})=\chi
_{1}^{l_{T}-1}\tilde{A}_{0}(\chi _{0},\chi _{1})$. Lemma \ref{Theo
IO 1ere carte} shows that
\begin{equation*}
\int\limits_{\mathbb{R}}\int\limits_{0}^{\infty }A_{0}(\chi
_{0},\chi_{1})e^{\frac{i}{h}\chi_{0}\chi_{1}^{2}}d\chi_{0}d\chi_{1}
=-\frac{1}{2}h^{d_{T}}\left\langle
\mathcal{F}(x_{-}^{d_{T}-1})(\chi _{0}),\tilde{A} _{0}(\chi
_{0},0)\right\rangle +\mathcal{O}(h^{d_{T}+\frac{1}{2}}).
\end{equation*}
With $\Lambda (\chi _{0})=\mathcal{F}(x_{-}^{d_{T}-1})(\chi
_{0})$, substituting the definition of $\tilde{A}_{0}$, gives
\begin{equation*}
\left\langle \Lambda (\chi _{0}),\tilde{A}_{0}(\chi
_{0},0)\right\rangle =2\pi \int e^{izr}\Lambda (\chi
_{0}(t,r,\theta ))a(t,r\theta ,0)dtdrd\theta dz.
\end{equation*}
But by construction $\chi _{1}=r,$ and since we have localized the
amplitude near $T$
\begin{equation}
\left\langle \Lambda (\chi _{0}),\tilde{A}_{0}(\chi
_{0},0)\right\rangle =-e^{i\frac{\pi }{2}(d_{T}-1)}\Gamma (d_{T})
\int \left\langle
(tQ_{t}(\theta)-i0)^{-d_{T}},a(t+T,0)\right\rangle d\theta dt.
\end{equation}
We can now use a result of \cite{GU}, also used in \cite{KhD} and
\cite {BPU}. Since the propagator $\mathrm{Exp}
(\frac{i}{h}tP_{h})$ is a FIO associated to the Lagrangian
manifold of the flow,
\begin{equation*}
\Lambda =\{(t,\tau ,x,\xi ,y,\eta )\text{ }/\text{ }(x,\xi )=\Phi
_{t}(y,\eta ),\text{ }\tau =-p(x,\xi )\},
\end{equation*}
it's principal symbol in the coordinates $(t,y,\eta )$ is given by
the half-density
\begin{equation*}
\mathrm{exp}(i\int\limits_{0}^{t}p^{1}(\Phi_s(y,\eta))ds)|dtdyd\eta
|^{\frac{1}{2}}.
\end{equation*}
The representation of the propagator with the kernel
\begin{equation*}
\frac{1}{(2\pi
h)^{n}}\int\limits_{\mathbb{R}^{n}}e^{\frac{i}{h}(S(t,x,\eta
)-\left\langle y,\eta \right\rangle )}(\alpha (t,x,\eta )+h\alpha
_{1}(t,x,\eta ,h))dy,
\end{equation*}
leads to the half-density
\begin{equation*}
\alpha (t,x,\eta )|dtdxd\eta |^{\frac{1}{2}}=\alpha (t,x,\eta )\frac{%
|dtdyd\eta |^{\frac{1}{2}}}{|\det (S_{x,\eta }^{\prime \prime })|^{\frac{1}{2%
}}}.
\end{equation*}
For our unique critical point $z_{0}$ we obtain
\begin{equation}
\alpha (t,z_{0})=|\det (S_{x,\eta }^{\prime \prime })(t,z_{0})|^{\frac{1}{2}%
}\exp (i\int\limits_{0}^{t}p^{1}(\Phi _{s}(z_{0}))ds)=\exp
(itp^{1}(z_{0})). \label{density}
\end{equation}
If there is no period of $d\Phi _{t}(z_{0})$ on
$\mathrm{supp}(\hat{\varphi})$ Theorem 5.6 of \cite{GU} gives
\begin{equation*}
\mathrm{Tr} (\psi ^{w}(x,hD_{x})\varphi (\frac{%
P_{h}-E_{c}}{h})\theta (P_{h})) \simeq \frac{1}{2\pi }e^{i\frac{\pi }{2}%
m_{0}}\int\limits_{\mathbb{R}}\frac{\hat{\varphi}(t)\exp (itp^{1}(z_{0}))}{%
|\det (\mathrm{Id}-d\Phi _{t}(z_{0})|^{\frac{1}{2}}}dt,
\end{equation*}
for a certain $m_0$. The denominator has a zero of order $d_{T}$
in $t=T$, hence for all $t\neq T$ in a sufficiently small
neighborhood of $T$, we have
\begin{equation*}
a(t,0)=\frac{(t-T)^{d_{T}}}{|\det (\mathrm{Id}-d\Phi _{t}(z_{0})|^{\frac{1}{2}}}%
\hat{\varphi}(t)\exp (itp^{1}(z_{0})),
\end{equation*}
Finally, since the contribution is smooth in $t$, we obtain that
\begin{gather*}
K(T) =-\frac{\Gamma (d_{T})}{2}\frac{\exp (i\frac{\pi
}{4}\mathrm{sgn} (q_{T}))}{(2\pi )^{1+d_{T}}| \det
(q_{T})|^{\frac{1}{2}}}
\exp (i\pi \frac{d_{T}-1}{2}\mathrm{sign}(Q_{T})), \\
\Lambda (\varphi ) =\left\langle
(t-T-i0)^{-d_{T}},\frac{(t-T)^{d_{T}}\exp (itp^{1}(z_{0}))}{|\det
(\mathrm{Id}-d\Phi
_{t}(z_{0})|^{\frac{1}{2}}}\hat{\varphi}(t)\right\rangle ,
\end{gather*}
and this completes the proof of Theorem 3.
\begin{remark}
\rm{If there are no rational relations between the eigenvalues of $Q^{+}$ and $%
Q^{-}$ all contributions of the non-zero periods of $d\Phi
_{t}(z_{0})$ are given by Theorem 3, as is shown by Proposition
\ref{Hess g, restric definie}. This gives a total contribution
\begin{equation*}
\sum\limits_{T\in \rm{supp}(\hat{\varphi})\backslash \{0\}}\gamma
(E_{c},T,h)=\sum\limits_{T\in \rm{supp}(\hat{\varphi})\backslash
\{0\}}(K(T)\Lambda _{T}(\varphi )+\mathcal{O}(h^{\frac{1}{2}})),
\end{equation*}
where the summation is over the non-zero periods of $d\Phi
_{t}(z_{0})$.}
\end{remark}
\textbf{The case $Q_{T}$ indefinite.}\newline Here we must study
carefully the solutions of the implicit function theorem. If
$z_{0}$ is zero then, by uniqueness, we have $z_{2}(t,0)=0$ for
all $t$, and
\begin{equation*}
\left\{
\begin{array}{l}
d_{z_{2}}\Psi (t,z_{1},z_{2}(t,z_{1}))=0, \\
\det (\mathrm{Hess}_{z_{2}}\Psi (T,0))\neq 0.
\end{array}
\right.
\end{equation*}
For our unique critical point $z_{1}=0$ we obtain
\begin{equation*}
\left\{
\begin{array}{c}
\dfrac{\partial z_{2}}{\partial t}(t,0)=-(\partial
_{z_{2},z_{2}}^{2}\Psi
(t,0,0))^{-1}\partial _{t}\partial _{z_{2}}\Psi (t,0), \\
\dfrac{\partial z_{2}}{\partial z_{1}}(t,0)=-(\partial
_{z_{2},z_{2}}^{2}\Psi (t,0,0))^{-1}\partial
_{z_{1},z_{2}}^{2}\Psi (t,0).
\end{array}
\right.
\end{equation*}
Equations (22) show that $\partial _{t}z_{2}(T,0)=0$ and $\partial
_{z_{1}}z_{2}(T,0)=0.$ Hence near $(T,0)$ we have
$z_{2}(t,z_{1})=\mathcal{O}(||(t-T,z_{1})||^{2})$. We use, again,
a decomposition w.r.t. $T$
\begin{equation}
\Psi (t,z_{1},z_{2}(t,z_{1}))=\Psi
(T,z_{1},z_{2}(T,z_{1}))+(t-T)g(t,z_{1},z_{2}(t,z_{1})).
\end{equation}
Since we have $d_{z_{2}}g(t,z_{1},z_{2}(t,z_{1}))=0,$ it follows
that
\begin{equation*}
\mathrm{Hess}_{z_{1}}(g(T,z_{1},z_{2}(T,z_{1}))_{|z_{1}=0}=Q_{T}.
\end{equation*}
Proposition \ref{TheopourRk} gives by identification of
homogeneous terms of same degree
\begin{equation}
\Psi
(T,z_{1},z_{2}(T,z_{1}))=R_{k}(z_{1},0)+\mathcal{O}(||z_{1}||^{k+1}).
\end{equation}
The stationary phase method in $z_{2}$, see formula
(\ref{Report8}), shows that
\begin{equation*}
I(T,h)\sim \sum\limits_{\nu=0}^{N}c_{\nu}h^{n+\nu-d_{T}}\int A_{\nu}(t,z_{1},0)e^{\frac{i%
}{h}(R(t,z_{1})+(t-T)g(t,z_{1}))}dtdz_{1},\text{ }h\rightarrow 0.
\end{equation*}
Copying the construction for a total period we obtain normal forms
for the phase $\Psi (t,z_{1})=R(t,z_{1})+(t-T)g(t,z_{1}) $, with
the decomposition w.r.t. $C_{Q_{T}}$ :
\begin{equation*}
A_{0}(t,z_{1},0)e^{\frac{i}{h}\Psi (t,z_{1})}\simeq \left\{
\begin{array}{l} A_{1,i}(\chi _{0},\chi _{1})e^{\pm
\frac{i}{h}\chi _{0}\chi _{1}^{2}}\text{
outside }C_{Q_{T}}, \\
A_{2,j}(\chi _{0},\chi _{1},\chi _{2})e^{\frac{i}{h}(\chi _{0}\chi
_{1}^{2}\chi _{2}\pm \chi _{1}^{k})}\text{ near
}C_{Q_{T}}\backslash
C_{R_{k}}, \\
A_{3,l}(\chi _{0},\chi _{1},\chi _{2},\chi
_{3})e^{\frac{i}{h}(\chi _{0}\chi _{1}^{2}\chi _{2}\pm \chi
_{1}^{k}\chi _{3})}\text{ near }C_{Q_{T}}\cap C_{R_{k}}.
\end{array}
\right.
\end{equation*}
Now the dimension is $l_{T}=2d_{T}$ and results obtained for total
periods show again that the contributions of normal forms 2 and 3
are dominating these of normal forms 1. Combining this with the
leading term of the stationary phase method gives
\begin{equation*}
I(T,h)=(2\pi h)^{n-d_{T}}\frac{e^{i\frac{\pi
}{4}\mathrm{sgn}(q_{T})}}{|\det
q_{T}|^{\frac{1}{2}}}(\int\limits_{\mathbb{R\times R}^{l_{T}}}
a(t,z_{1},0)e^{\frac{i}{h}(R(t,z_{1})+(t-T)g(t,z_{1}))}dtdz_{1}+\mathcal{O}(h)),
\end{equation*}
and the contribution of a non-total period is computed by
restriction of all objects to $\frak{F}_{T}$. This proves Theorems
\ref{Main result Q positive} and \ref{Main result Q non positive}
in their general forms.
\section{Examples}
\noindent \textbf{Perturbation of harmonic oscillators.} Let be
\begin{equation}
H(x,\xi )=\frac{1}{2}((x_{1}^{2}+\xi _{1}^{2})-(x_{2}^{2}+\xi
_{2}^{2}))+(x_{2}^{2}+\xi _{2}^{2})^{2},
\end{equation}
with critical energies $0$ and $-\frac{1}{16}$. We consider
$\Sigma _{0}=\{(x,\xi )$ $/$ $H(x,\xi )=0\}$, the origin is the
only singularity of this surface. The system is integrable with
\begin{equation*}
\Phi _{t}(x,\xi )=\left(
\begin{array}{c}
x_{1}\cos (t)+\xi _{1}\sin (t) \\
x_{2}\cos ((4(x_{2}^{2}+\xi _{2}^{2})-1)t)+\xi _{2}\sin
((4(x_{2}^{2}+\xi
_{2}^{2})-1)t) \\
\xi _{1}\cos (t)-x_{1}\sin (t) \\
\xi _{2}\cos ((4(x_{2}^{2}+\xi _{2}^{2})-1)t)-x_{2}\sin
((4(x_{2}^{2}+\xi _{2}^{2})-1)t)
\end{array}
\right) .
\end{equation*}
On $\Sigma _{0}$ the function $4(x_{2}^{2}+\xi _{2}^{2})-1$ is
bounded and some elementary arithmetical considerations show that
we can find a neighborhood $V(T)$ of the origin, in $\Sigma _{0},$
such that all periodic trajectory has a period greater than $|T|$.
Hence $(H_{5})$ is true here. We have
\begin{equation*}
d\Phi _{t}(0)=\left(
\begin{array}{cccc}
\cos (t) & 0 & \sin (t) & 0 \\
0 & \cos (t) & 0 & -\sin (t) \\
-\sin (t) & 0 & \cos (t) & 0 \\
0 & \sin (t) & 0 & \cos (t)
\end{array}
\right) ,
\end{equation*}
and $2\pi$ is a total period. Expanding the flow in a
Taylor-series gives
\begin{equation*}
d_{0}^{3}\Phi _{2\pi }((x,\xi )^{3})=\left(
\begin{array}{c}
0 \\
4\xi _{2}(x_{2}^{2}+\xi _{2}^{2}) \\
0 \\
-4x_{2}(x_{2}^{2}+\xi _{2}^{2})
\end{array}
\right) .
\end{equation*}
Hence $R_{4}(x,\xi )$ is negative on $C_{Q}\backslash \{0\}$,
since
\begin{equation*}
R_{4}(x,\xi )=\frac{1}{24}<(x,\xi ),Jd_{0}^{3}\Phi _{2\pi }((x,\xi )^{3})>=-%
\frac{1}{6}(x_{2}^{2}+\xi _{2}^{2})^{2}.
\end{equation*}
\newline \textbf{Siegel and Moser example.} To illustrate some properties on periodic
trajectories, Siegel \& Moser introduced in \cite{S-M} the
Hamiltonian
\begin{equation}
H(x,\xi )=\frac{1}{2}(x_{1}^{2}+\xi _{1}^{2})-(x_{2}^{2}+\xi
_{2}^{2})+x_{1}\xi _{1}x_{2}+\frac{1}{2}(x_{1}^{2}-\xi
_{1}^{2})\xi _{2}.
\end{equation}
The origin is the only critical point on $\Sigma _{0}=\{(x,\xi )$ $/$ $%
H(x,\xi )=0\}$ and trajectories from the surface $\{x_{1}=\xi _{1}=0\}$ are $%
2\pi $-periodic$.$ With $p=(x_{1}^{2}+\xi _{1}^{2}),$
$q=(x_{2}^{2}+\xi _{2}^{2}),$ we obtain, $p^{\prime \prime
}=4pq+p^{2}$, (see \cite{S-M}). Since $p$ is strictly positive,
$p$ is strictly convex and hence non-periodic. Consequently, on
$\Sigma _{0}$ there is no non-trivial periodic trajectories. The
linearized flow at the critical point is
\begin{equation*}
d\Phi _{t}(0)=\left(
\begin{array}{cccc}
\cos (t) & 0 & \sin (t) & 0 \\
0 & \cos (2t) & 0 & \sin (2t) \\
-\sin (t) & 0 & \cos (t) & 0 \\
0 & -\sin (2t) & 0 & \cos (2t)
\end{array}
\right),
\end{equation*}
with a resonance of order 3 since, $2w_{1}-w_{2}=0$. For the
''period'' $2\pi $ we have
\begin{equation*}
R_{3}(x,\xi )=\frac{1}{6}\left\langle -Jd_{0}^{2}\Phi _{2\pi
}((x,\xi )^{2}),(x,\xi )\right\rangle =\frac{\pi
}{4}(2x_{1}x_{2}\xi _{1}+\xi _{2}(x_{1}^{2}-\xi _{1}^{2})).
\end{equation*}
The intersection $C_{Q}\cap C_{R_{3}}$ is obtained by solving the
system
\begin{equation*}
\left\{
\begin{array}{c}
(x_{1}^{2}+\xi _{1}^{2})-2(x_{2}^{2}+\xi _{2}^{2})=0, \\
2x_{1}x_{2}\xi _{1}+\xi _{2}(x_{1}^{2}-\xi _{1}^{2})=0.
\end{array}
\right.
\end{equation*}
This leads to the surfaces
\begin{equation*}
(S_{1}):\left\{
\begin{array}{c}
x_{2}(x_{1},\xi_{1})=\dfrac{x_{1}^{2}-\xi_{1}^{2}}{\sqrt{2}x_{1}
\xi_{1}}\sqrt{\dfrac{x_{1}^{2}\xi _{1}^{2}}{x_{1}^{2}+\xi _{1}^{2}}} \\
\xi _{2}(x_{1},\xi _{1})=-\sqrt{2}\sqrt{\dfrac{x_{1}^{2}\xi
_{1}^{2}}{ x_{1}^{2}+\xi _{1}^{2}}}
\end{array}
\right. ,\text{ }(S_{2}):\left\{
\begin{array}{c}
x_{2}(x_{1},\xi _{1})=-\dfrac{x_{1}^{2}-\xi _{1}^{2}}{\sqrt{2}x_{1}\xi _{1}}%
\sqrt{\dfrac{x_{1}^{2}\xi _{1}^{2}}{x_{1}^{2}+\xi _{1}^{2}}} \\
\xi _{2}(x_{1},\xi _{1})=\sqrt{2}\sqrt{\dfrac{x_{1}^{2}\xi _{1}^{2}}{%
x_{1}^{2}+\xi _{1}^{2}}}
\end{array}
\right. .
\end{equation*}
By symmetry we just examine gradients on the first surface
$(S_{1})$, we have
\begin{equation*}
\nabla R_{3|S_{1}}=\frac{\pi }{4}\left(
\begin{array}{c}
-\sqrt{2}x_{1}\xi _{1}^{2}\sqrt{\dfrac{x_{1}^{2}+\xi
_{1}^{2}}{x_{1}^{2}\xi
_{1}^{2}}} \\
2x_{1}\xi _{1} \\
\sqrt{2}x_{1}^{2}\xi _{1}\sqrt{\dfrac{x_{1}^{2}+\xi
_{1}^{2}}{x_{1}^{2}\xi
_{1}^{2}}} \\
(x_{1}^{2}-\xi _{1}^{2})
\end{array}
\right) ,\text{ }\nabla Q_{|S_{1}}=\left(
\begin{array}{c}
2x_{1} \\
4\dfrac{x_{1}^{2}-\xi _{1}^{2}}{\sqrt{2}x_{1}\xi_{1}}\sqrt{\dfrac{
x_{1}^{2}\xi _{1}^{2}}{x_{1}^{2}+\xi _{1}^{2}}} \\
2\xi _{1} \\
-4\sqrt{2}\sqrt{\dfrac{x_{1}^{2}\xi _{1}^{2}}{x_{1}^{2}+\xi
_{1}^{2}}}
\end{array}
\right) .
\end{equation*}
But, since the minor determinant extracted from $\nabla
Q_{|S_{1}}$ and $\nabla R_{3|S_{1}}$
\begin{equation*}
D(x_{1},\xi _{1})=\left| \left(
\begin{array}{cc}
x_{1} & \xi _{1} \\
-x_{1}\xi _{1}^{2}\sqrt{\dfrac{x_{1}^{2}+\xi _{1}^{2}}{x_{1}^{2}\xi _{1}^{2}}%
} & x_{1}^{2}\xi _{1}\sqrt{\dfrac{x_{1}^{2}+\xi
_{1}^{2}}{x_{1}^{2}\xi _{1}^{2}}}
\end{array}
\right) \right| =\sqrt{\frac{x_{1}^{2}+\xi _{1}^{2}}{x_{1}^{2}\xi _{1}^{2}}}%
x_{1}\xi _{1}(x_{1}^{2}+\xi _{1}^{2}),
\end{equation*}
is non zero for $(x_{1},\xi _{1})\neq 0$, $\nabla Q$, $\nabla
R_{3}$ are linearly independent on $C_{Q}\cap C_{R_{3}}\cap
\mathbb{S}^3$.


\begin{thebibliography}{9}
\bibitem{BB} R.Balian and C.Bloch, Solution of the Schr\"odinger equation in term
of classical path, Annals of Physics \textbf{85} (1974) 514-545.
\bibitem{BPU} R.Brummelhuis, T.Paul and A.Uribe, Spectral estimate near a critical level,
Duke Math. Journal \textbf{78} (1995) no. 3, 477-530.
\bibitem{BU} R.Brummelhuis and A.Uribe, A semi-classical trace formula for Schr\"odinger
operators, Commun. Math. Phys. \textbf{136} (1991) no. 3, 567-584.
\bibitem{C-P} A.M.Charbonnel and G.Popov, A semi-classical trace
formula for several commuting operators, Commun. Partial
Differential Equations \textbf{24} (1999) no. 1-2, 283-323.
\bibitem{DUI1} J.J.Duistermaat, Oscillatory integrals Lagrange immersions and unfolding
of singularities, Commun Pure Appl. Math. \textbf{27} (1974),
207-281.
\bibitem{GU} V.Guillemin and A.Uribe, Circular symmetry and the trace formula,
Inven. Math. \textbf{96} (1989) no. 2, 385-423.
\bibitem{GUT} M.Gutzwiller, Periodic orbits and classical quantization conditions,
J. Math. Phys. \textbf{12} (1971).
\bibitem{HOR1} L.H{\"o}rmander, "The analysis of linear partial operators 1,2,3,4",
Springer-Verlag (1985).
\bibitem{HOR2} L.H{\"o}rmander, "Seminar on singularities of solutions of linear
partial differential equations", Annals of mathematical studies
91, Princeton University Press (1979) 3-49 .
\bibitem{KhD} D.Khuat-Duy, A semi-classical trace formula for Schr\"odinger operators in the
case of a critical energy level, Th\`ese de l'universit\'e Paris
9.
\bibitem{KhD1} D.Khuat-Duy, A semi-classical trace formula for Schr\"odinger operators in the
case of a critical energy level, J. Funct. Anal. \textbf{146}
(1997) no. 2, 299-351.
\bibitem{S-M} J.K.Moser and C.L.Siegel, "Lectures on celestial mechanics", Springer-Verlag (1995).
\bibitem{PU} T.Paul and A.Uribe, Sur la formule semi-classique des
traces, Comptes Rendus S\'eances Acad.Sci. S\'erie I. \textbf{313}
(1991) no. 5, 217-222.
\bibitem{P-P} V.Petkov and G.Popov, Semi-classical trace formula and
clustering of the eigenvalues for Schr\"odinger operators, Ann.
Inst. Henri Poincar\'e Phys. Th\'eorique \textbf{68} (1998) no. 1,
17-83.
\bibitem{Rob} D.Robert,"Autour de l'approximation
semi-classique", Progress in mathematics Volume 68, Birkh\"{a}user
Boston (1987).
\bibitem{WON} R.Wong, "Asymptotic approximations of integrals",
Academic Press Inc. (1989).
\end{thebibliography}
\end{document}